\documentclass[lettersize,journal]{IEEEtran}
\usepackage{amsmath,amsfonts}
\usepackage{algorithm}
\usepackage{algorithmic}
\usepackage{array}
\usepackage{bm}
\usepackage{booktabs}
\usepackage{graphicx}
\usepackage[caption=false,font=normalsize,labelfont=sf,textfont=sf]{subfig}
\usepackage{textcomp}
\usepackage{stfloats}
\usepackage{url}
\usepackage{verbatim}
\usepackage{enumitem}
\usepackage{graphicx}
\usepackage{multicol}
\def\BibTeX{{\rm B\kern-.05em{\sc i\kern-.025em b}\kern-.08em
    T\kern-.1667em\lower.7ex\hbox{E}\kern-.125emX}}
\usepackage{balance}
\usepackage{tikz-cd}
\usepackage[capitalise]{cleveref}

\crefname{remark}{Remark}{Remarks}

\newcommand{\vd}{{\bm d}}
\newcommand{\vs}{{\bm s}}
\newcommand{\vn}{{\bm n}}
\newcommand{\vp}{{\bm p}}
\newcommand{\vr}{{\bm r}}
\newcommand{\vv}{{\bm v}}
\newcommand{\vu}{{\bm u}}
\newcommand{\vx}{{\bm x}}
\newcommand{\vy}{{\bm y}}
\newcommand{\vz}{{\bm z}}

\newcommand{\vS}{{\bm S}}
\newcommand{\vJ}{{\bm J}}
\newcommand{\vU}{{\bm U}}

\newcommand{\vV}{{\bm V}}
\newcommand{\Id}{{\bm I_d}}

\newcommand{\vmu}{{\bm \mu}}

\newtheorem{remark}{Remark}

\newcommand{\mua}{{\bm{\delta\mu}}_a}

\def\lp{\left(}
\def\rp{\right)}

\def\la{\left\langle}
\def\ra{\right\rangle}

\DeclareMathOperator*{\argmin}{argmin}
\DeclareMathOperator*{\soft}{soft}
\DeclareMathOperator*{\sign}{sign}

\begin{document}
\title{A Learned--SVD approach  for  Regularization  in Diffuse Optical Tomography}
\author{
Alessandro Benfenati
\thanks{Alessandro Benfenati is with Department of Environmental and Science Policy, Via Celoria 2, {Milan}, {20133}, {Italy}, and with Gruppo Nazionale Calcolo Scientifico, INDAM}
Giuseppe Bisazza 
\thanks{Giuseppe Bisazza is with Department of Mathematics, Via Saldini 50, {Milan}, {20133}, {Italy}}
Paola Causin 
\thanks{Paola Causin is with Department of Mathematics, Via Saldini 50, {Milan}, {20133}, {Italy}, and with Gruppo Nazionale Calcolo Scientifico, INDAM}}

\markboth{Transaction on Computational Imaging}{Alessandro Benfenati, Giuseppe Bisazza,\MakeLowercase{\textit{(et al.)}: A Learned SVD approach for Inverse Problem Regularization in Diffuse Optical Tomography}}

\maketitle

\begin{abstract}
	Diffuse Optical Tomography (DOT) is an emerging technology in medical imaging which 
	employs  light in the NIR spectrum to estimate the distribution of optical coefficients in biological tissues
	for diagnostic and monitoring purposes.
	DOT reconstruction implies the solution of a severely ill--posed inverse problem,
	for which regularization techniques are mandatory in order to achieve reasonable results. 
	Traditionally, regularization techniques put a variance prior on the desired solution/gradient via  regularization parameters, whose choice requires a fine tuning, specific for each case.   
	In this work we explore deep learning techniques in a fully data--driven approach, able 
	of reconstructing the generating signal (optical absorption coefficient) in an automated way. 
	We base our approach on the so-called Learned Singular Value Decomposition, which has been proposed for general inverse problems, and we tailor it to the DOT application. 
	We perform tests with increasing levels of noise on the measure, and compare it with standard variational approaches.
\end{abstract}

\begin{IEEEkeywords}
	Diffuse Optical Tomography, Inverse Problems, SVD, Deep Learning, Regularization 
\end{IEEEkeywords}


\section{Introduction}
\label{sec:intro}
\noindent\IEEEPARstart{D}{iffuse} Optical Tomography (DOT) 
is an emerging technology which employs light in the NIR spectral window to investigate biological tissues
for diagnostic and monitoring purposes. 
Its use has been explored in different medical fields, including imaging of brain, thyroid, prostate, 
breast cancer screening and, more in general,  
for detecting and monitoring body functional changes related to blood flow~\cite{hoshi2016overview,yamada2014diffuse,jiang2018diffuse}. 
Tomographic reconstruction in DOT aims at recovering the spatial distribution of tissue 
optical properties of the investigated organ, which can be related 
to chromophores concentration 
(mainly water, oxy and deoxy-hemoglobin and, in the case of adipose tissue as breast, lipids) for functional evaluation.
 To do this, light emitted
from external sources at different positions is let to propagate throughout 
the investigated organ and the emerging photon flux is measured
on the tissue boundary~\cite{arridge1999optical}. DOT reconstruction is a notoriously severely ill-posed and ill-conditioned
inverse problem since at NIR wavelengths the outgoing light consists of a mixture of very few  coherent and quasi--coherent photons and a predominant component of incoherent (diffusive) photons
which experienced multiple scattering events.  As such, the DOT technology 
is intrinsically sensitive to measurement noise and model errors.

Model-based image reconstruction algorithms (also known as variational approaches) 
have represented in the past the standard approach to
perform DOT reconstruction (inverse problem). They usually consist of three components~\cite{hoshi2016overview}:\begin{enumerate}[label=\roman*)]
	\item a forward mathematical model that provides a prediction of the
	measurements based on a guess of the system parameters
	(the optical coefficients: the absorption coefficient $\mu_a$ [cm$^{-1}$] and the scattering coefficient $\mu_s$
	[cm$^{-1}$]); 
	\item an objective functional that compares the predicted data with the measured data;
	\item an efficient way of
	updating the system parameters of the forward model,
	which in turn provides a new set of predicted data. 
\end{enumerate}

This latter point has been approached both by linearization via the Born or Rytov approximations~\cite{Boas01} and by nonlinear iterations formulated as optimization problems~\cite{Arridge2009optical}. 
Due to the ill-conditioned nature of the reconstruction problem, all these approaches include
a form of regularization, which provides hard or soft prior information on the solution. 
As for this latter category, the $\ell_2$--norm (Tikhonov/ridge regression) penalization is the benchmark classic regularization approach and has been widely used in DOT literature~\cite{5728925,doi:10.1137/090781590,10.1117/12.2230074}. 
The $\ell_1$-norm (lasso) penalization has been mainly used in this context to improve
sharp edges detection, obtained as a by--product 
of sparsity enhancement~\cite{Okawa:11}. 
Some of the authors of the present paper proposed in~\cite{Causin19,Causin20} the use of 
an Elastic--Net regularization term.  
This approach was investigated both in 2D and 3D domains under the Rytov approximation and provided 
results superior to the use of pure $\ell_2$  or $\ell_1$ strategies,  the combination  
of which being able to give a satisfactorily stable solution with a preserved sparsity pattern. 
An alternative approach was proposed in~\cite{Benfenati2020}, based on 
Bregman iterations~\cite{Zhang201120}: this iterative technique consists in substituting the regularization function with its Bregman distance from the previous iterate. It has been proved \cite{Benfenati14} that this procedure provides a solution to the optimization problem and moreover it possess a remarkable contrast enhancement ability \cite{Benfenati6,Benfenati11}, confirmed also in DOT framework \cite{Benfenati2020}.
Other regularization functionals have also been explored, among which we cite here the combined use of a TV method
with $\ell_1$--norm in the {\tt TOAST++} software for DOT reconstruction~\cite{Schweiger2014toast++}.
All these approaches are far from being flawless: poor reconstruction results 
are common outcomes already in mildly complex situations and/or the reconstruction 
requires a very heavy computational time. 
Specifically, a main critical point is the fine tuning of regularization parameters required
to obtain a reasonable solution. These parameters were observed to be strongly 
dependent on the geometry, mesh-size, optical coefficients distribution and level 
of noise of each single case, so that a general, effective, strategy cannot be obtained.

Recently, the outstanding
performance  on computer vision tasks of deep learning (DL) algorithms based on Neural Networks (NNs), and especially Convolutionary Neural Networks (CNNs), 
has motivated studies to explore their use also in DOT reconstruction,
with the aim to overcome the above discussed weaknesses. The application of DL techniques 
to DOT problems is still at its beginnings and very recent: 
in~\cite{Wu20} a \emph{Regularization by Denoising} (RED) approach was proposed, while \cite{Yoo2020} trained a NN for inverting in an end--to-end fashion the Lippman--Scwhinger equation, by firstly learning the pseudoinverse of the nonlinear mathematical operator modelling the DOT physics,  and then applying an encoder--autoencoder network to remove the 
remaining artifacts. In~\cite{mozumder2021model}, the authors combined DL with a model-based approach using Gauss--Newton iterations where the update function was learned via a CNN.  We refer to~\cite{zhang2019brief} for a review of the 
other few applications of DL in DOT reconstruction.  
In this work, we investigate and adapt to DOT reconstruction the general strategy  for inverse problems of~\cite{Boink19}, 
called Learned--SVD. This approach consists in mapping measured data and generating signal 
into feature spaces and bridge them via a ``singular value operator`` which acts as a learned regularization,
and alleviates from the burden of~\emph{a--priori} choosing the regularization. 
The idea of learning the regularization via NNs is not new in the general field of inverse problem
solution: for example, in~\cite{Lunz18} a regularization term was learned using only unsupervised data;  
in~\cite{Kobler_2020_CVPR, Kobler20} a Total Deep Variation functional was learned, while the NETT framework  of~\cite{Li2020} employs as regularizer the encoder part of an encoder--decoder NN. However, 
at the best of our knowledge, none of these approaches -including the present Learned-SVD- has been applied to DOT reconstruction.  

We examined the performance of the Learned-SVD and compared it to
classic variational approaches. According to the measured metrics, 
the present NN-based strategy performs significantly better and increasing 
levels of noise on the data have a reduced impact than on 
variational methods. 

This paper is organized as follows. \Cref{sec:setting} presents the abstract mathematical setting
of the DOT reconstruction problem; \Cref{sec:variational} presents the classic variational approaches,
with specific focus on the Rytov linearization strategy; \Cref{sec:DL} presents the Learned-SVD strategy with its connection to the SVD decomposition and Tikhonov regularization; \Cref{sec:numexp} presents the results obtained with the proposed approach together with a comparison with standard variational reconstruction approaches; \Cref{sec:concl} draws the conclusions. 

\section{General  setting}
\label{sec:setting}
The goal of the DOT reconstruction is to map 
measured light data collected on the investigated tissue boundary into accurate approximate spatial  distribution of the optical
coefficients inside the tissue itself. Due to relevant aspects, which directly impact on the mathematical reconstruction problem, 
in this paper we focus on DOT based on continuous waves (CW) systems, which
are among the most widely used in clinical settings for breast cancer screening, object of our past and
present research work. In these systems, which are typically embedded in low-to-mid cost range hardware, 
the light source continuously emits light into the tissue at a single frequency and the 
measured data are light fluence outgoing from the tissue boundary.
A minimal CW--DOT setting consists thus in an array of $n_s$ pointwise light sources 
located at positions $\vs$ and emitting light at NIR wavelength into the tissue and a set of $n_d$ measures of the outgoing light flux.  
These latter can be collected by detectors physically located on the tissue surface at positions
$\vd$ or by a charge--coupled device (CCD) camera which produces (indirect) measures of the outgoing
light intensity via grey-levels in an image. 
Additional elaborations are usually performed on the raw measures in order to (at least partially) filter out noise, eliminate
under- and over-saturated pixels by thresholding procedures and define an appropriate ROI. 
We do not enter into details of these procedure, referring to~\cite{vavadi2016automated} 
for an overview of the implicated workflow. 

Let $q \in Q$ denote the set of spatially dependent optical properties
of the tissue, $u \in V$ the intensity (or fluence) of the light in the tissue, and let $f=f({\vs}) \in W$ be a given source
located at position ${\vs}$, where  
$Q, V$  and $W$ are the Banach spaces of all admissible optical properties, light field and source terms, 
respectively. 
The set of parameters $q$ along with the given sources $f$ determine observations of quantities related
to the variable state $u$ that are collected at the detector locations $\vd$ as  
\begin{equation}
y=\mathcal{M}(q,u; f),
\label{eq:ymes}
\end{equation}
where $y=[y(\vd_1,\vs_1),\dots,y(\vd_{n_d},\vs_{n_s})], $ is the measurement vector 
and where  $\mathcal{M}: (Q \times V \times W) \rightarrow Y$ is a bounded linear functional,  $Y$ being a Banach 
data space. 
The operator $\mathcal{M}$ is a general representation of the measurement.  
Notice that measured data are invariably affected by noise, mainly electronic noise due to
the detector acquisition chain but also from other artifacts like (in realistic conditions) patient's 
movements and ambient light. 

\medskip

On the other side, we introduce the parameter-to-state
map  $\mathcal{F}: \Theta \rightarrow Y$ such that
\begin{equation}
\mathcal{F}(\theta) = 
\mathcal{P}u.
\label{eq:paramtostate}
\end{equation}
where $\theta$ is a set of parameters that will be specified below. 
One can pursue an approximation of~$\mathcal{F}$ via a traditional (physics-driven, see Sect.~\Cref{sec:variational})
or machine learning--based (data--driven, see Sect.~\Cref{sec:DL}) approach.
In the first case, $\Theta \equiv Q$ and 
the parameter--to--state map arises from an underlying mathematical model,
which, in the DOT context, describes light propagation in a biological tissue. 
Such a model can be represented in an abstract manner as   
\begin{equation}
\mathcal{A}(q,u)=f,
\label{eq:forward_abs}
\end{equation}
where $\mathcal{A}: (Q\times V) \rightarrow W$ and if $u$ is a solution to~\eqref{eq:forward_abs} for a given $q$ and given source $f=f(\vs)$ then $\mathcal{P}u$ is the trace of the solution and it is a quantity directly comparable to  the datum $y$  in~\eqref{eq:ymes}.
On the other hand, when considering data--driven approaches, one does not dispose of an explicit 
model as in~\eqref{eq:forward_abs}, but directly learns $\mathcal{F}$ from data, that is  
learning the set of parameters~$\Theta$ (weights of the net).
In both cases, to close the problem one must introduce a non--negative
function, called discrepancy or loss function, $\mathcal{D}: \Theta \times Y \rightarrow \mathbb{R}^+$  
that - considering the case of least squares - can be written as 
\begin{equation}
\mathcal{D}(\theta;y)=\frac12 \Vert \mathcal{F}(\theta) 
- y\Vert_{Y}^2
\end{equation}
and the DOT reconstruction problem is then formulated as
 \begin{equation}
\begin{array}{l}
{\theta}^*=\displaystyle \argmin_{\theta  \in \Theta}\mathcal{D}(\theta; y).
\end{array}
\label{eq:DOT_abs}
\end{equation}

\section{Variational approaches for DOT reconstruction}
\label{sec:variational}
\subsection{Mathematical models of light propagation}

As anticipated in the previous section, classic variational approaches in DOT reconstruction rely on the existence 
of a mathematical model.  In this case, one disposes of an operator 
$\mathcal{A}(\cdot,\cdot)$ in~\eqref{eq:forward_abs} 
which is typically a differential law 
 augmented by appropriate boundary conditions  describing 
 light propagation in scattering media.
A physically reasonable and cost--effective model of this kind 
is derived by performing an expansion of the 
Radiative Transfer Equation (RTE) in spherical harmonics (see {\em e.g.}~\cite{Durduran2010DiffuseOF} for details and derivation). 
The following diffusion equation (DE) model is then
obtained in steady state:\\
\begin{equation}\label{eq:eq_diff}
 -\nabla\cdot (D(\vr) \nabla u(\vr))
+ \mu_a(\vr) u(\vr) 
 = f(\vs),
\end{equation}
where $\vr$ is the position [cm] in the connected open domain~$\Omega\subset \mathbb{R}^d, d=2$ or $3$, 
representing the tissue to be investigated (for example the breast), with outward unit normal vector ${\vn}$.
In~\eqref{eq:eq_diff},  $u=u(\vr)$ represents photon fluence [W/cm$^{2}$]
due to the light source $f=S_0\delta(\vr-\vs)$  [W/cm$^{3}$],  
$S_0$ being the source intensity. 
The parameter space is represented by:
\begin{equation}\label{param:RTE:Q}
Q = \left\{ q=(\mu_a,\mu_s)\,:\, 0 < \mu_a \le \overline{\mu}_a,
 0 < \mu_s \le \overline{\mu}_s\right\},
\end{equation}
$\overline{\mu}_a$ and $\overline{\mu}_s$ being positive upper bounds for the respective optical coefficient fields. 
When considering the CW technology is customary to assume $\mu_s=\overline{\mu}_s$ to be a known constant, so that also the diffusion coefficient $D$ [cm], defined as
\begin{equation}
 \label{eq:diff_coef}
D(\vr) = \frac{1}{3(\mu_a(\vr)+(1-g)\mu_s)} \approx \frac{1}{3 (1-g) \mu_s}, 
\end{equation}
is a constant,  with  $g$  anisotropic scattering factor depending on the specific tissue.
Model~\eqref{eq:eq_diff} is equipped with boundary conditions on the domain boundary $\partial \Omega$
\begin{equation}
\label{eq:BCs}
    D\nabla u \cdot {\vn} +\frac{1}{2A_c} u  =0,  \\
\end{equation}
which describe the behavior of light at the interface of the 
tissue with air or another material according to the Fresnel law,  
$A_c$ being the accomodation coefficient at the tissue--air interface (see, {\em e.g.},~\cite{Mansuripur2002}).

\subsection{DOT reconstruction via Rytov linearization}

Rytov linearization is a common algorithm implemented in the reconstruction 
software embedded in DOT instrumentation, since it  
computes the solution in a short time with  limited request for processing 
power.  We follow here 
the derivation provided by Ishimaru~\cite[Vol. II, Ch. 17]{Ishimaru78}, with a  slight modification 
to include the presence of volumetric sources.
In brief, 
we assume the linearization 
$\mu_a(\vr)=\mu_{a,0}(\vr) + \delta \mu_a(\vr)$, where $0 < \mu_{a,0} \le \overline{\mu}_a $ is the 
background value and 
$\delta \mu_a$ a perturbation term corresponding to  
localized contrast regions in the tissue. 
We also introduce the exponential change of variables  
\begin{equation}
u(\vr)=u_0(\vr){e^{\psi_1(\vr)}},
\label{eq:expchange}
\end{equation}
where  $u_0(\vr)$ is the light fluence field in background conditions and $\psi_1=\log(u/u_0)$ is the 
{\em logarithmic amplitude fluctuation} of the light due to the presence of contrast regions.
Inserting~\eqref{eq:expchange} into~\eqref{eq:eq_diff},
performing algebraic manipulations (see~\cite{causin2020inverse} for a detailed derivation) 
and neglecting higher order terms, 
the following equation for  the combined quantity $(u_0 \psi_1)$ is obtained
\begin{equation}
\begin{array}{ll}
\displaystyle\left[\Delta - \frac{\mu_{a,0}}{D}\right](u_0 \psi _1)(\vr)=\frac{\delta \mu_a(\vr)}{D} u_0(\vr)& \quad {\rm in} \,\,\Omega,
\\\\
\displaystyle D\nabla (u_0 \psi_1) \cdot {\vn} +\frac{1}{2A_c} (u_0 \psi_1)  =0 & \quad {\rm on}\,\, \partial \Omega. \\
\end{array}
\label{eq:MH}
\end{equation} 

The DOT inverse problem is obtained
by solving for the optical coefficient distribution which minimizes the
discrepancy between the modeled and measured~$\vy$. 
To do this, the parameter-to-state map is obtained from~\eqref{eq:MH},
whose solution provides the field $\psi_1$ for source 
$\displaystyle\left(\frac{\delta \mu_a}{D} u_0\right)$, upon evaluation of
at the detectors locations through.  
As for the data, in the spirit of the Rytov perturbative approach,
one must dispose of a set of fluence measures ${u}_{0}$ which correspond to ``unperturbed/background conditions''
and a set of fluence measures $u$ which correspond to a condition for which the lightpaths express
the optic footprint of possible internal contrast regions. 
Notice that the definition of the perturbed status
depends on the specific technological setting: for example, in breast screening 
it can be obtained by compression of the breast.
This causes increased accumulation of deoxy-hemoglobin and in the possible
presence of tumor--related neo-vasculature, a differential increase is expected with 
locally altered absorption coefficient.   

\subsection{Numerical approximation}

The DOT problem is inherently discrete, since in practice one disposes only of a finite set of measures
of light in the detectors, obtained in correspondence of each of the $n_s$ sources, 
for a total of $M=n_d \times n_s$ measures. 
Moreover, problem~\eqref{eq:MH} is solved via numerical or semi-analytical techniques: 
finite element or finite volume methods can be
used (see, {\em e.g.}, \cite{arridge1999optical}) or alternatively, as in our implementation, 
meshless approaches based on Green's function method. 
In this latter case,  one performs the convolution between
the source term of~\eqref{eq:MH} and the Green’s function kernel $G(\vr-\vr^\prime),  
\vr, \vr^\prime \in \Omega$, 
for the modified Helmholtz
operator $[\Delta - \alpha^2]$, with attenuation coefficient
$\alpha =\sqrt{\mu_{a,0}/D}$.
Upon partitioning the domain into $N$ voxels, the following  sensitivity matrix 
is obtained 
\begin{equation}
\label{eq:sensit}
\vJ=J_{ij}=\left[\frac{\Delta V_j}{U_0(\vs_k,\vd_l)}
 G(\vd_l-\vr_j)\frac{1}{D} U_0(\vs_k,\vr_j)\right],
\end{equation}
where the index $i=1,\dots, M \rightarrow \{l,k\}, l=1, \dots, n_d, k=1, \dots, n_s$ stands
for the detector/source pair
and the index $j$ for the voxel number, $j=1,\dots,N$,
with $\Delta V_j$ volume  and $\vr_j$ centroid of the $j$-th voxel, respectively.  
Notice that~\eqref{eq:sensit} already corresponds to the evaluation
of the solution at the detectors locations and 
that the quantity $U_0(\cdot,\cdot)$ appearing in it is an approximation of
the background field evaluated at position as in its first argument for source in position as 
in its second argument and obtained as well from the Green's function
upon application of the reciprocity principle.  
One ends up thus with  
the following discrete parameter-to-state map 
\begin{equation}
\label{eq:sensit0}
\vJ \mua = {\boldsymbol \psi}_1,
\end{equation}
where $\mua \in \mathbb{R}^{N \times 1}$ is the vector of 
unknown perturbations of the absorption coefficient evaluated at
 $N$ locations (centroids) in the tissue volume discretized in voxels,
 ${\boldsymbol \psi}_1 \in \mathbb{R}^{M \times 1}$ is the vector of 
the fluctuations evaluated at the detectors' locations.
The solution of the DOT inverse problem is then obtained by 
the discrete minimization problem
\begin{equation}
\label{eq:inv}
\mua^* =\argmin_{\mua} \mathcal{D}(\vJ \mua;\vy),
\end{equation}
where  ${\vy}=\log{\left({u}/{u_0}\right)\vert_\vd} \in \mathbb{R}^{M \times 1}$ is the vector 
of measured fluctuations at the detectors locations for the different sources: it is the discrete measurements equivalent to $\boldsymbol\psi_1$ from \eqref{eq:expchange}.
Due to the severe ill-conditioning  of~\eqref{eq:inv}, regularization techniques are mandatory in order to obtain a physically coherent solution~\cite{Sciacca21,Panagiotou:09}. 
One is then led to consider the following regularized discrete variational formulation
\begin{equation}
\mua^{*} =\argmin_{\mua } \mathcal{D}(\vJ \mua;\vy) +
\alpha \mathcal{R}(\vmu_a),
\label{eq:rego}
\end{equation}
where $\mathcal{R}$ is a regularization function and  the regularization parameter $\alpha$ balances its influence on the solution. Both the choice of the regularization function and the regularization parameter are crucial.  A classical strategy is to choose $\mathcal{R}$ as the $\ell_2$--norm of the solution, which sets a uniform variance prior (known also as Tikhonov regularization or ridge regression in statistical frameworks), and to hand--pick $\alpha$ from a parametric analysis (\emph{e.g.}, GCV, L--curve). 
Choosing $\mathcal{R}$ as the $\ell_1$--norm of the solution promotes instead sparseness \cite{Cao:07}, preserving 
only a (relatively) small number of non--zero components. An interesting approach, which compromises
the benefits of the $\ell_2$ and $\ell_1$ penalizations is the Elastic-Net, where one sets
\begin{equation}
 \mathcal{R}(\vmu_a)=
\label{eq:regoEN}
\alpha(\vartheta||\vmu_a||_1 + (1-\vartheta)||\vmu_a||_2^2),
\end{equation}
$ \theta \in [0,1]$ being a parameter which weights the contribution of the two norms. 
By varying $\theta$, sparsity
is traded for the inclusion of more solution components contributions. Another performant approach is the Bregman technique, which consists in an iterative procedure where $\mathcal{R}$ is substituted with its Bregman distance at each iteration. More details on these two latter approaches are given in \cref{ssec:varappr}.

In general, the choice of $\alpha$ is quite
delicate. For pure $\ell_1$--norm sparsity-promoting penalization, there is no golden rule for an automatic
selection: some specific strategies refer to particular instances~\cite{meinshausen2006high}. For the pure $\ell_2$--norm
penalization, there exist more consolidated methods to select the parameter: a
general possibility, implemented in advanced software packages, is
to use $k$--fold cross validation techniques. 

\medskip

\begin{remark} The Rytov approach yields a linear instance
of the parameter-to-state  mapping $\mathcal{F}$ from optical parameters to light fluence values
based on a physical model.  
Nonlinear operators can also be obtained from other models/approximations. For example, a more accurate formulation could be derived from the  Lippman--Schwinger equation, which amounts to consider also higher order terms in the Rytov perturbative approach (see~\cite{yodh1995spectroscopy,Yoo2020}). Alternatively, one may want to consider the full RTE equations,
albeit at a significantly higher computational cost. 
Problem~\eqref{eq:inv} can thus be written in the more general form (notice that here we 
consider the full absorption coefficient and not only the perturbation and $\vy$ represents in general
the measure): 
\begin{equation}
\vmu_a^{*} =\argmin_{\vmu_a \in \mathbb{R}^{N \times 1}} \mathcal{D}(\mathcal{F}(\vmu_a);\vy). 
\label{eq:inv1}
\end{equation}
One should observe that regularization is in any case necessary and can be expressed analogously 
to~\eqref{eq:rego}.
\end{remark}

\section{Neural--network based approaches for DOT reconstruction}
\label{sec:DL}
\begin{figure*}[!t]
\begin{centering}
	\subfloat[\label{fig:structures_BOINK}]{\includegraphics[width=.70\textwidth]{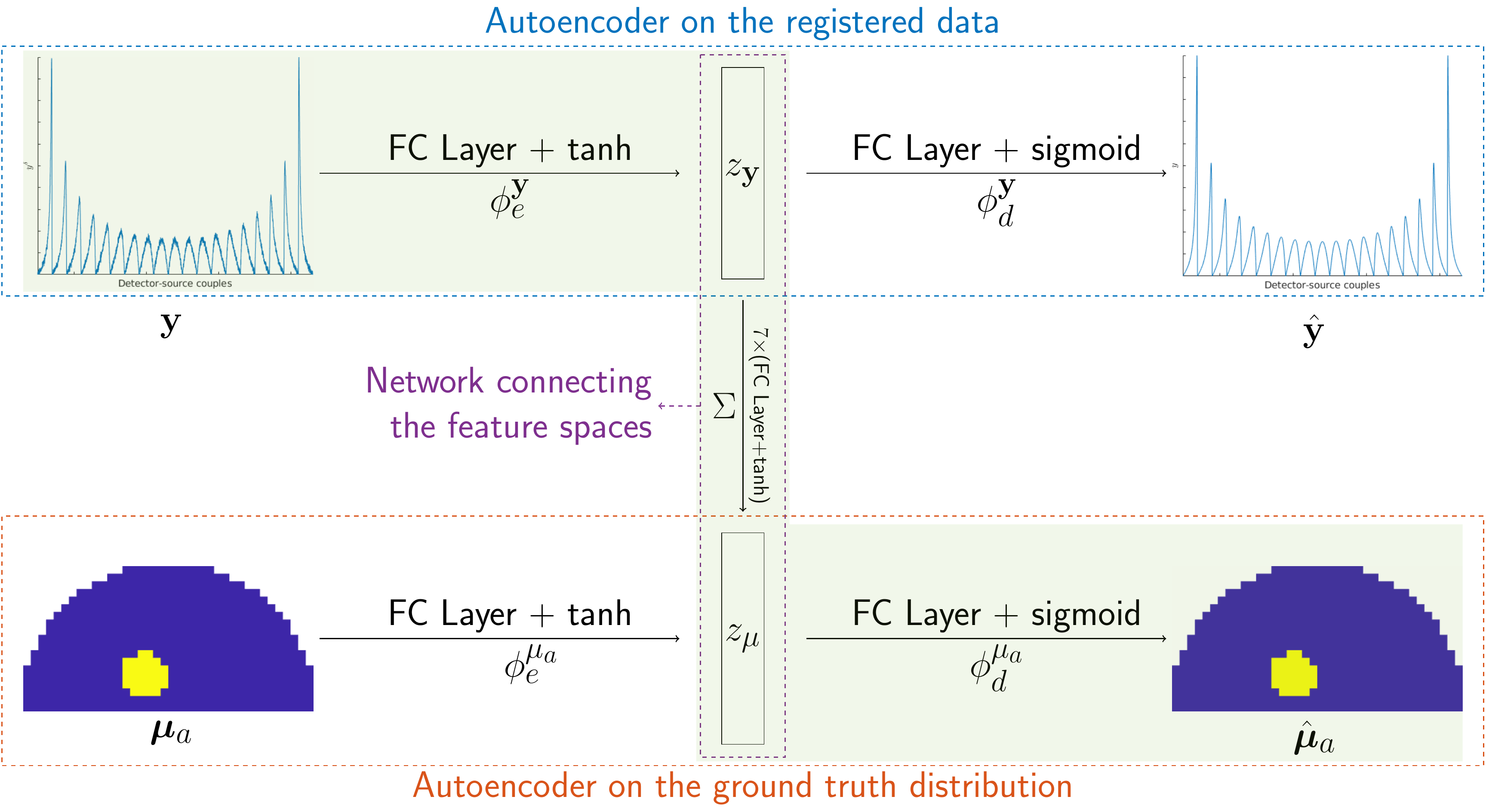}}\hfill\subfloat[\label{fig:workflow}]{\includegraphics[width=0.22\textwidth]{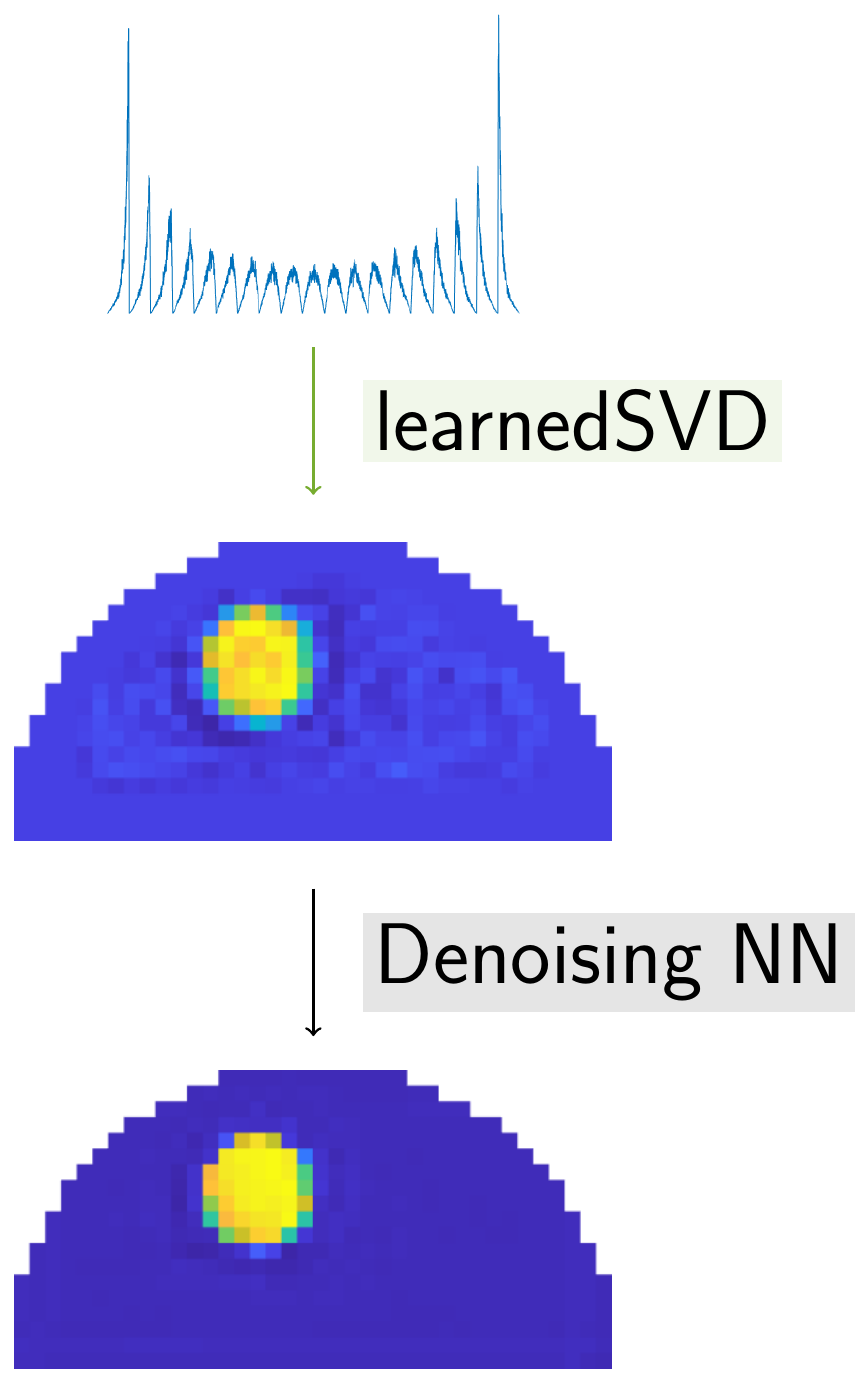}}
\end{centering}
\caption{(a): Learned-SVD strategy for DOT reconstruction. The upper part of the scheme represents the dAE
which receives as an input data $\vy$ and outputs the reconstruction $\hat{\vy}$. In doing this, denoise
is also performed.  
Data are ordered in detector-source pairs as specified in the caption of Fig.~\ref{fig:domain}. 
The bottom part represents the sAE which reconstructs the noise--free ground truth of the originating signal. The operator $\Sigma$ is learned in order to connect the two feature representations. The end-to-end workflow is enlightened in green and defines the mapping $\phi_d^\mu\circ\Sigma\circ\phi_e^y$; (b): test samples flow through the trained Learned-SVD NN and are further denoised 
to improve the results as in~\cite{Yoo2020}. Abbreviations: FC  = Fully Connected layer, tanh = hyperbolic tangent activation function, sigmoid = sigmoid activation function.}
\label{fig:networks}
\end{figure*}
We follow here the Learned-SVD strategy introduced in~\cite{Boink19}, which aims
at mimicking via a NN the SVD decomposition  when its truncated or regularized version is employed. 
In order to understand the concept at the basis of this approach, it is useful to draw a parallel
with the models introduced in the previous section. 
In this case, when
the parameter-to--state operator $\mathcal{F}$ is linear, as it is for example using the Rytov strategy,
one may write $\mathcal{F}=\vU{\vS}\vV^\top$, where ${\vS}$ is a diagonal matrix containing the singular values and write the solution by the inversion formula $\mathcal{F}^{-1}=\vV {\vS}^{-1} \vU^\top$. 
This idea is extended  here to more general nonlinear parameter-to-state operators,
which are not formally known from a model but directly learnt by data.   
Namely, one  learns a generalized SVD 
decomposition via a special class of NNs called autoencoder (AE) architectures.  
To do this, similarly to~\cite{Boink19}, we consider two AEs: 
the first, called data-AE (dAE), acts on data $\vy$ 
and is defined by 
\begin{equation}
\vz_\vy = \phi_e^y(\vy), \qquad  \hat{\vy} = \phi_d^y(\vz_\vy), 
\end{equation}
where $\vz_\vy$ is the encoded feature (latent representation) of the noisy data $\vy$ and $\hat{\vy}$ its decoded version. In the SVD framework, $\phi_e^\vy$ and $\phi_d^\vy$ can be seen as the correspective of $\vU^\top$ and $\vU$ for the linear case, respectively. This encoder--decoder has also by construction a denoising behaviour since the 
AE is trained in order to yield a noise--free reconstruction of the noisy input. The second AE, called signal--AE (sAE), acts on the generating signal and is defined by 
\begin{equation}
\vz_{\mu} = \phi_e^\mu({\vmu_a}), \qquad \widehat{\vmu}_a = \phi_d^{\mu}({\vmu_a}),
\end{equation}
where $\vz_\mu$ is the encoded feature (latent representation) of the generating signal $\vmu_a$ and $\widehat{\vmu}_a$ 
its decoded version. Within the SVD framework, the decoder $\phi_d^\mu$ can be considered as the counterpart of $\vV$. The two latent feature representations $\vz_\vy$ and $\vz_{\mu}$ are related by a bridging operator $\Sigma$, such that $\vz_{\mu}=\Sigma(\vz_\vy)$, which plays the scaling role of the singular values in~$\vS$ in the classical SVD approach. 

The strategy consists thus in: 
\begin{enumerate}[label=\roman*)]
\item Encoding the data $\vy$ in the latent space via the encoder $\phi_e^y$ (corresponds to perform 
the product $\vU^\top\vy$).
\item Connecting the latent variables $\vz_\vy$ and $\vz_\mu$ via the operator $\Sigma$ (mimics the computation $\vS^{-1}\vU^\top\vy$).
\item Decoding the latent variables $\vz_\mu$ into the coefficient distribution via the decoder $\phi_d^\mu$
(corresponds to the final left multiplication by $\vV$ in the SVD).  
\end{enumerate} 
 The parallelism of the proposed scheme with classical SVD is further represented by
 the scheme in \eqref{eq:lSVDstruct}.
\begin{equation}
\label{eq:lSVDstruct}
\begin{tikzcd}
\vy \arrow{r}{\phi_e^y}[swap]{\vU^\top} &  \vz_y \arrow{d}{\Sigma}[swap]{\vS}\arrow{r}{\phi_d^y}[swap]{\vU}   &\hat{\vy}\\
\vmu_a \arrow{r}{\phi_e^\mu}[swap]{\vV^\top} &  \vz_\mu\arrow{r}{\phi_d^\mu}[swap]{\vV} &\widehat{\vmu}_a\\
\end{tikzcd}
\end{equation}
Fig.~\ref{fig:networks} depicts the application of the  Learned-SVD to the DOT context, which  
results in a fully {\em data--driven} method.

\medskip

\begin{remark}
	\label{rem:tikSVD}
Tikhonov regularization and (truncated) SVD share a common role in regularization. Indeed, consider (\ref{eq:rego}) with a linear
operator and when $\mathcal{R}$ is chosen as $\frac12\|\cdot\|_2^2$: the solution to this problem can be expressed using SVD as
$$
\mua^* = \sum_i  f_\alpha(\sigma_i)\la \vu_i, \vy \ra\,\vv_i,
$$
where $\vv_i$ and $\vu_i$ are the $i$--th right and left eigenvectors and $\sigma_i$ is the $i$--th singular value of $\vJ$, respectively,
and where the function $f_\alpha(\sigma)=\frac{\sigma}{\sigma^2+\alpha}$  filters out the smaller singular values,
according to the threshold $\alpha$ (regularization parameter).
While in the classical approach this parameter is chosen a-priori and used in the specific function 
$f_\alpha(\sigma)$, the Learned-SVD approach learns the operator $\Sigma$ which plays the role of the filtering function $f$ and doses automatically 
the amount of required regularization. In doing so, it learns at the same time the regularization functional $\mathcal{R}$ and the parameter $\alpha$, which are both encompassed in $\Sigma$.
\end{remark}

\medskip

\begin{remark}
Notice that in the Learned-SVD approach, $\vy$ directly represents the measurements $u\vert_\vd$ and not the logarithmic amplitude fluctuation, since no linearization has been performed. As a consequence, the reconstruction yields
the full absorption coefficient vector. 
\end{remark}

\section{Results}
\label{sec:numexp}
We assessed the performance of the proposed NN-based approach via quantitative indices and we compared it  
to the results obtained with variational regularization obtained via the Elastic-Net 
penalization or the Bregman approach.   
Experiments were carried out on a synthetic phantom both in noise--free and noisy conditions. 
All the numerical experiments were carried on a laptop PC with operating system Linux~22.04, endowed with an Intel(R) Core(TM) i5–8250U CPU (1.60GHz) and 16 GiB RAM memory. The NN architectures were implemented in MatLab\textregistered R2022a
using the Deep Learning Toolbox\textsuperscript{TM}.

\subsection{Generation of synthetic data}
\begin{figure}[!t]
	\newcommand{\factor}{0.4}
	\begin{center}
	\includegraphics[width=\factor\textwidth]{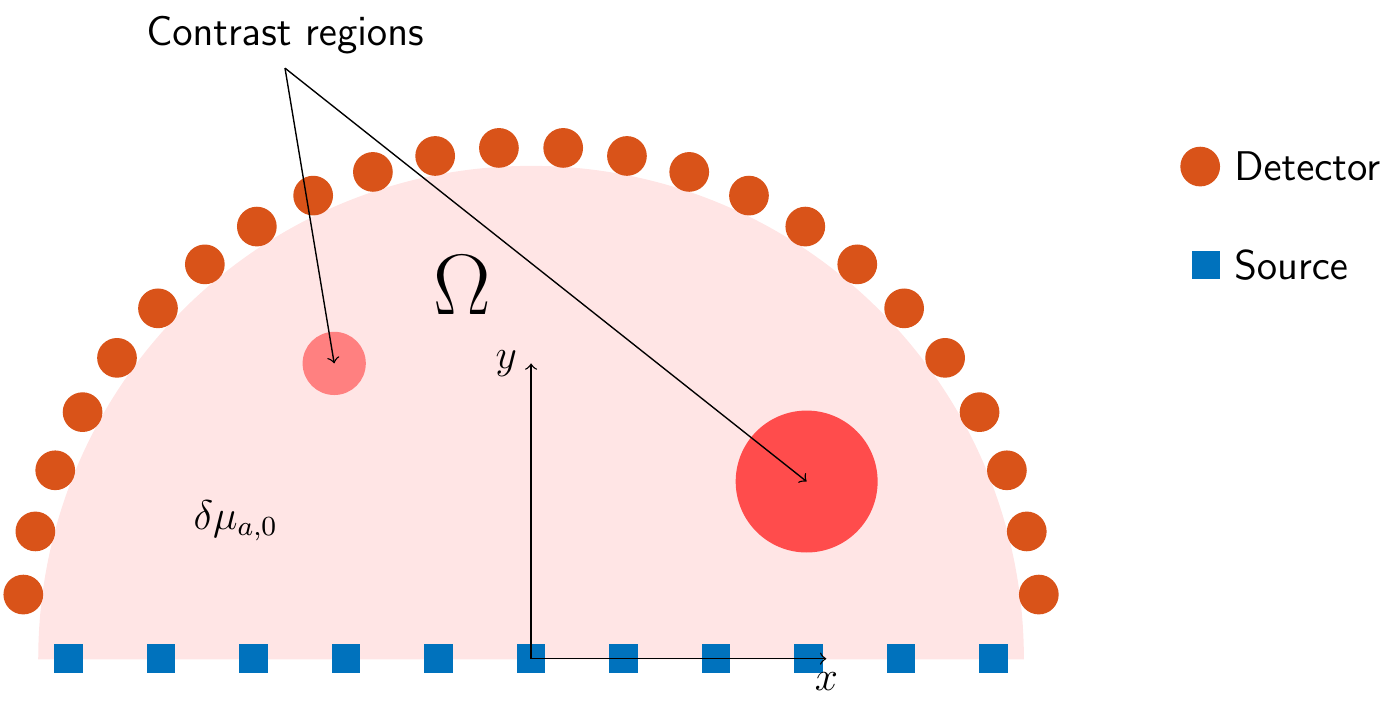}
	\end{center}
	\caption{Geometry of the domain $\Omega$ considered in the numerical simulations. 
	The detectors are uniformly disposed
	along the curved boundary with numeration progressing in counterclockwise sense 
	starting from the rightmost detector; the sources are aligned along the bottom boundary and
	positioned 1mm inside~$\Omega$, with  numeration increasing from left to right.
	The absorption 
	coefficient is $\mu_{a,0}$ except for the contrast regions which have
	increased absorption coefficient. }
\label{fig:domain}
\end{figure}

\noindent  We considered  synthetic measures  generated  via the DE
forward model discretized by an in--house finite element code using
a fine mesh. 
The domain consists in a 2D semidisk with radius 5~cm, 
with 19 light sources positioned 1~mm inside the straight boundary and 20 detectors uniformly distributed on the semicircular portion of the boundary. The background coefficient  was set to $\mu_{a,0}=0.01$~cm$^{-1}$ and 
we considered constant $(1-g)\mu_s=0.1$~cm$^{-1}$ (see Fig.~\ref{fig:domain})
For reconstruction, we discretized the domain into square voxels with side 0.25~cm. 
To train the network, we generated a set of 1500 samples, each including one or two circular contrast regions with 
random radius and position inside the domain. The 
contrast regions were
chosen to have absorption coefficient equal to $3,4$ or $5$-fold the background coefficient. When there were two contrast regions,  each region was allowed to have different
absorption coefficient.
The test dataset used to evaluate the performance of the network and to compare it with 
the variational approaches contained 150 samples built according to the same rules as above.  
In our experiments, we considered noise-free data as well as noisy measures, obtained by adding white Gaussian noise to the fluence values at the detectors with variance 1\%, 3\% or 5\%, respectively. 
\begin{figure*}[!t]
\newcommand{\factor}{0.18}
\begin{center}
	\subfloat[\label{fig:gt}]{\includegraphics[width=\factor\textwidth]{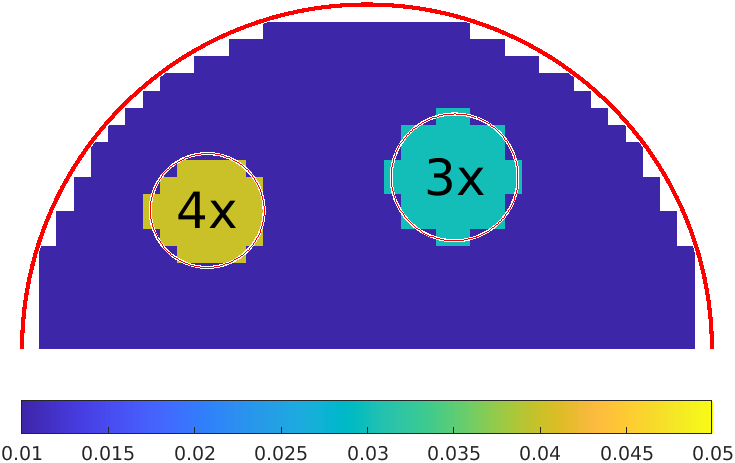}\hspace{0.025\textwidth}\includegraphics[width=\factor\textwidth]{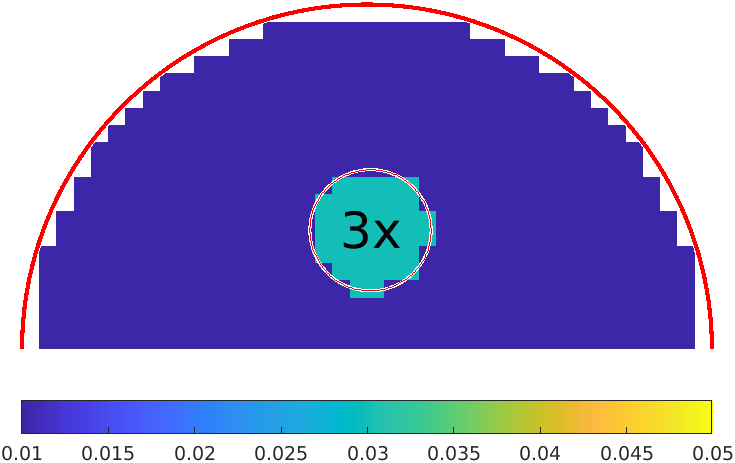}\hspace{0.025\textwidth}\includegraphics[width=\factor\textwidth]{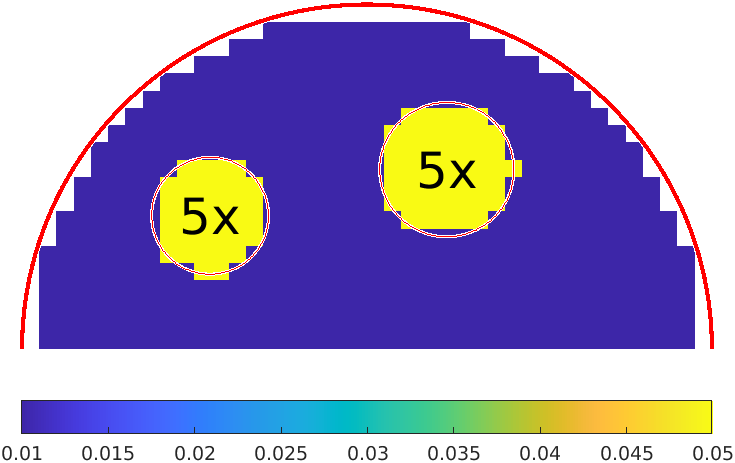}\hspace{0.025\textwidth}\includegraphics[width=\factor\textwidth]{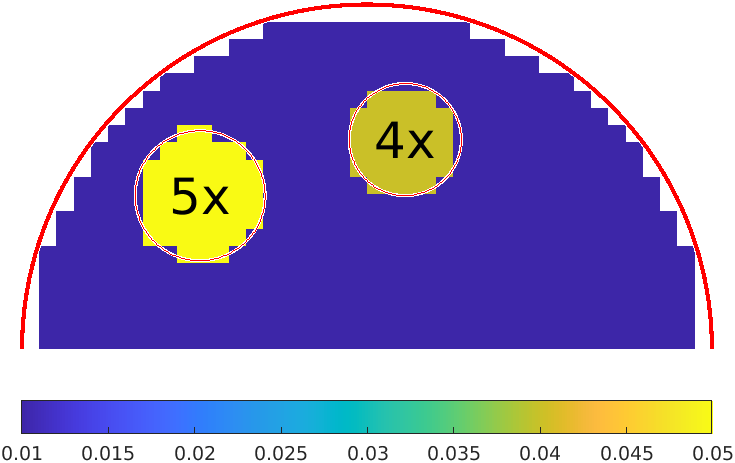}\hspace{0.025\textwidth}\includegraphics[width=\factor\textwidth]{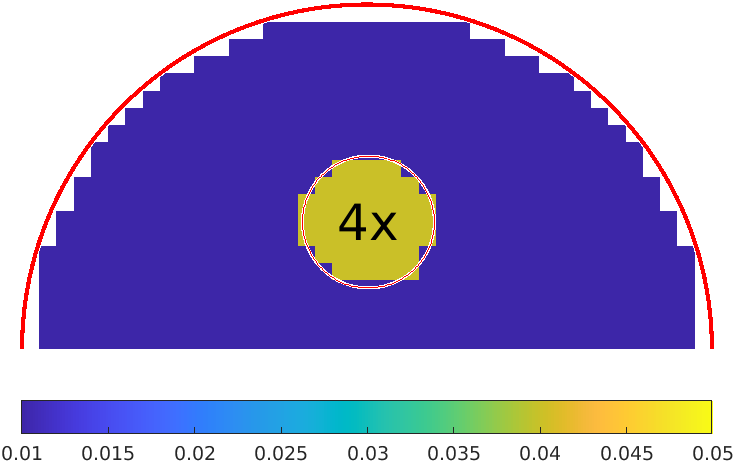}}
	
	\subfloat[\label{fig:nfree}]{\includegraphics[width=\factor\textwidth]{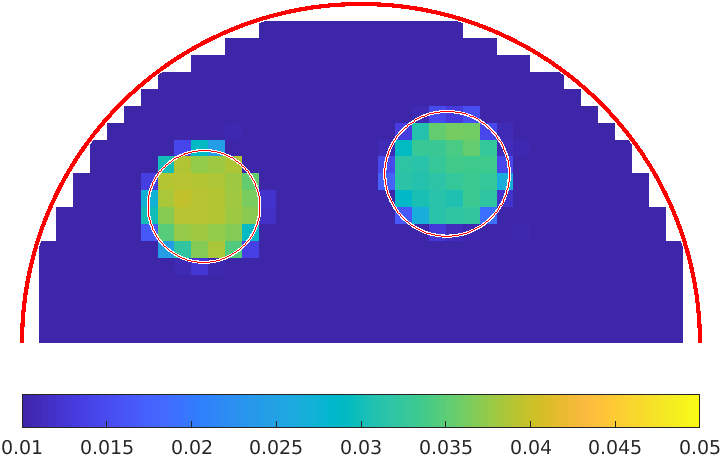}\hspace{0.025\textwidth}\includegraphics[width=\factor\textwidth]{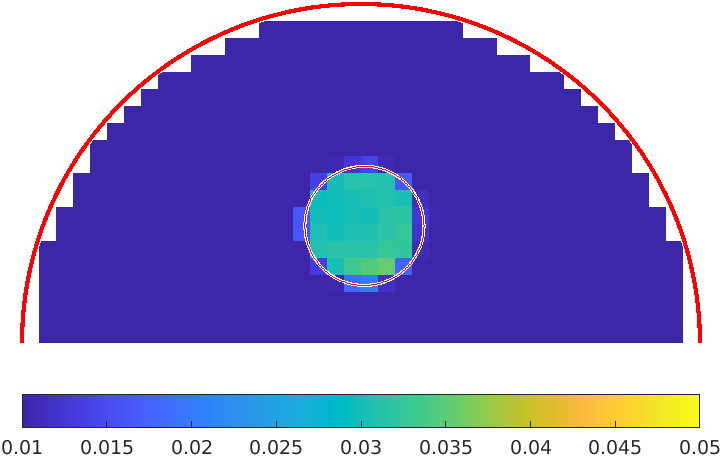}\hspace{0.025\textwidth}\includegraphics[width=\factor\textwidth]{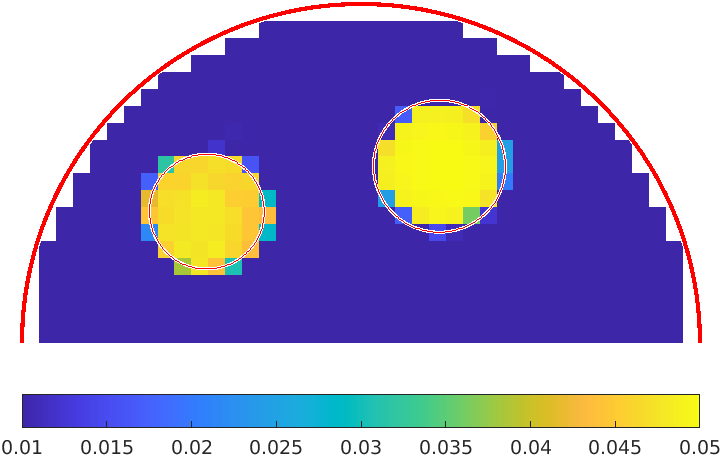}\hspace{0.025\textwidth}\includegraphics[width=\factor\textwidth]{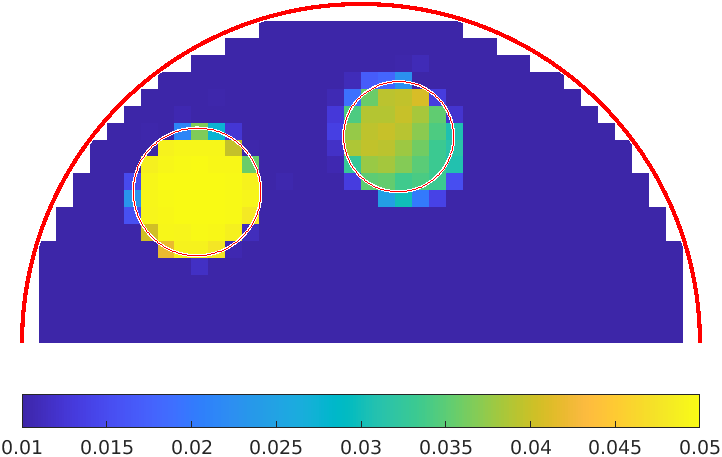}\hspace{0.025\textwidth}\includegraphics[width=\factor\textwidth]{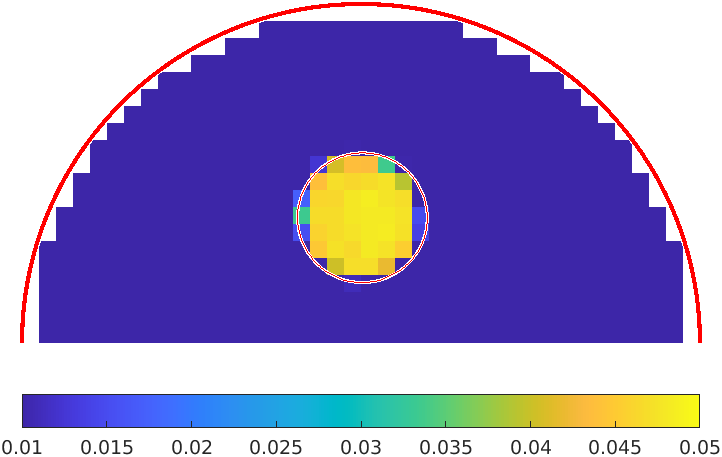}}
	
	\subfloat[\label{fig:G1}]{\includegraphics[width=\factor\textwidth]{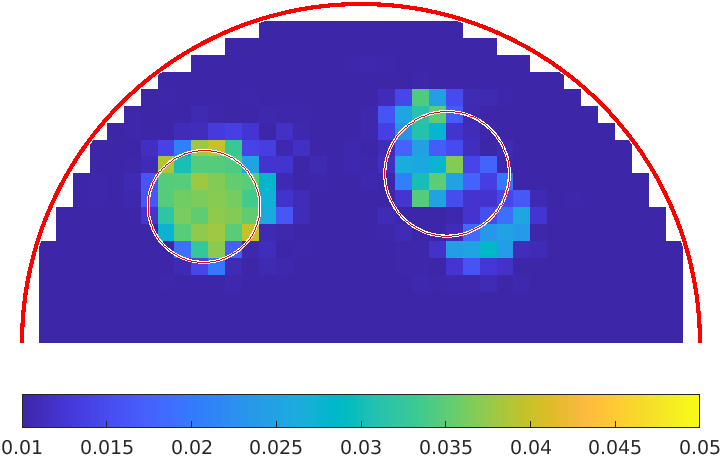}\hspace{0.025\textwidth}\includegraphics[width=\factor\textwidth]{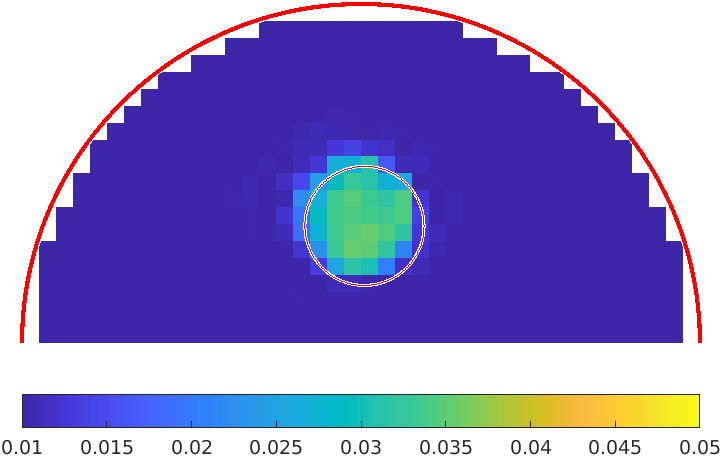}\hspace{0.025\textwidth}\includegraphics[width=\factor\textwidth]{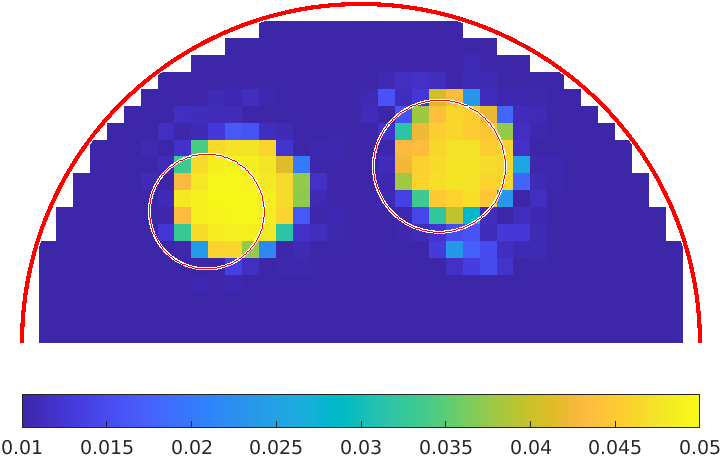}\hspace{0.025\textwidth}\includegraphics[width=\factor\textwidth]{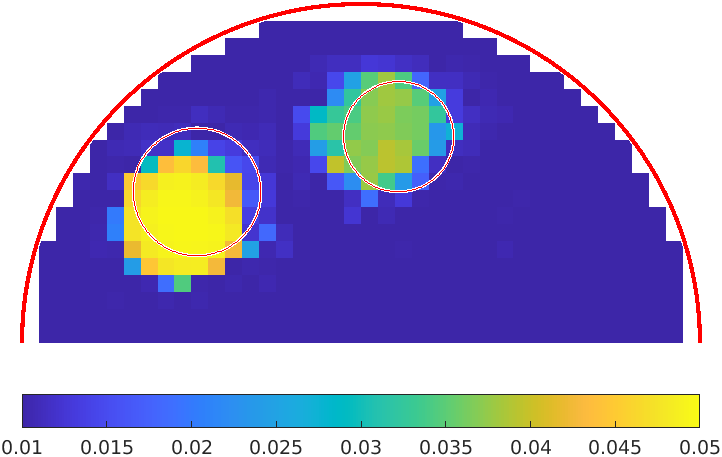}\hspace{0.025\textwidth}\includegraphics[width=\factor\textwidth]{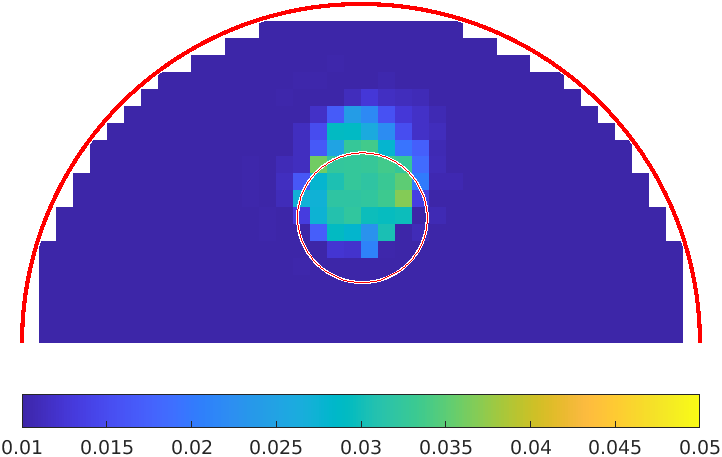}}
	
	\subfloat[\label{fig:G3}]{\includegraphics[width=\factor\textwidth]{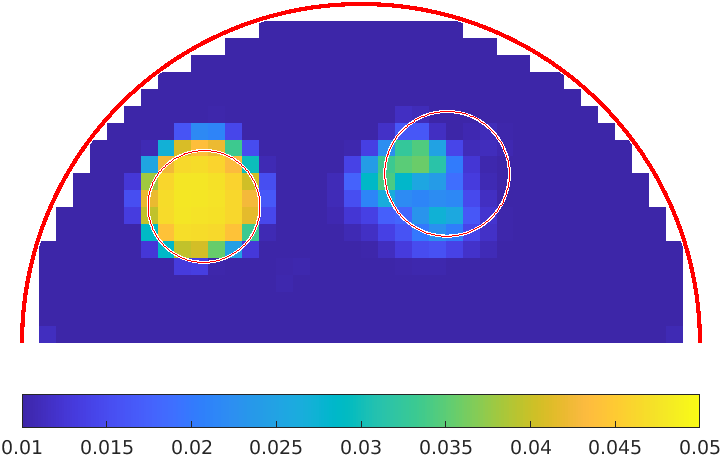}\hspace{0.025\textwidth}\includegraphics[width=\factor\textwidth]{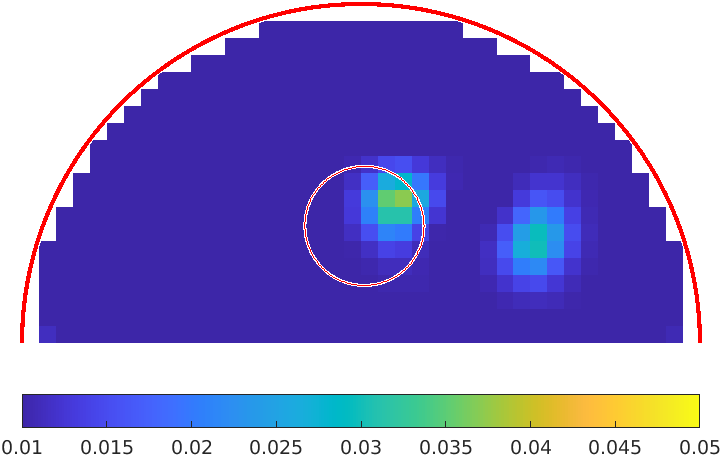}\hspace{0.025\textwidth}\includegraphics[width=\factor\textwidth]{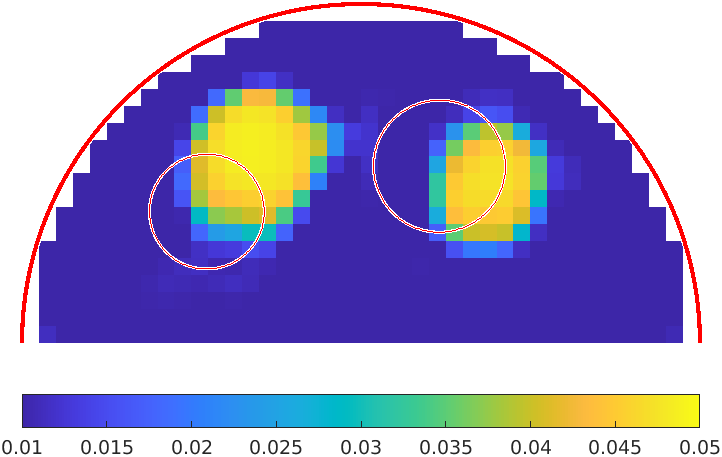}\hspace{0.025\textwidth}\includegraphics[width=\factor\textwidth]{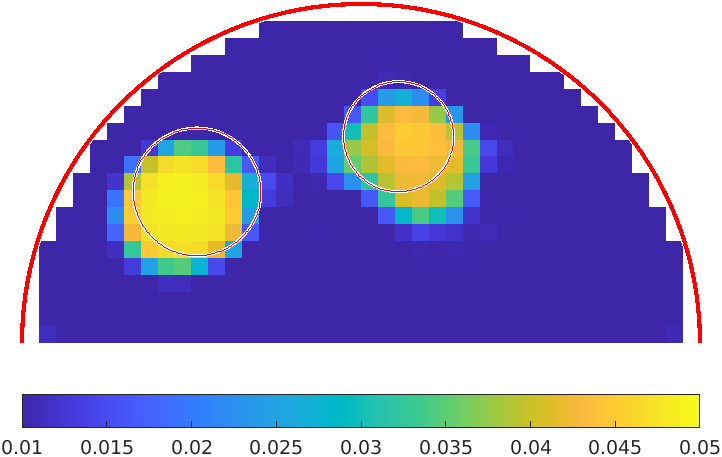}\hspace{0.025\textwidth}\includegraphics[width=\factor\textwidth]{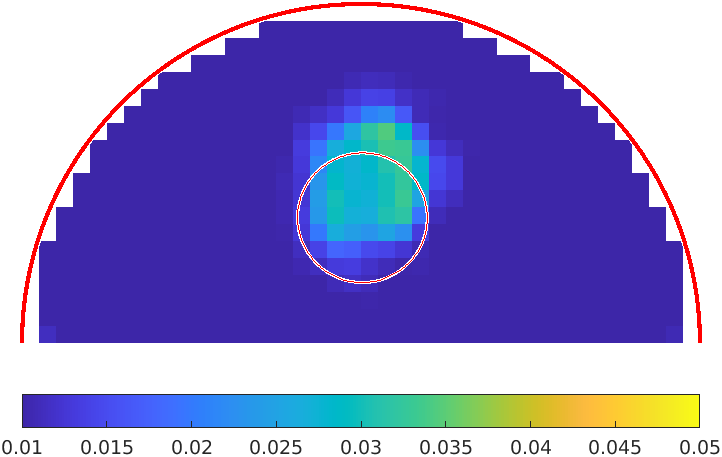}}
	
	\subfloat[\label{fig:G5}]{\includegraphics[width=\factor\textwidth]{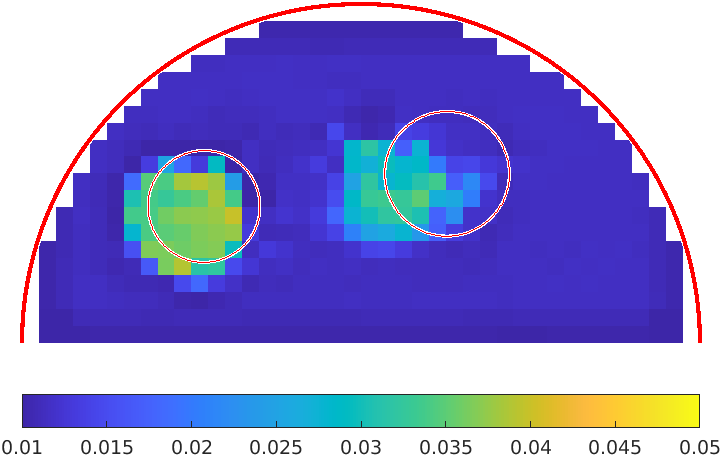}\hspace{0.025\textwidth}\includegraphics[width=\factor\textwidth]{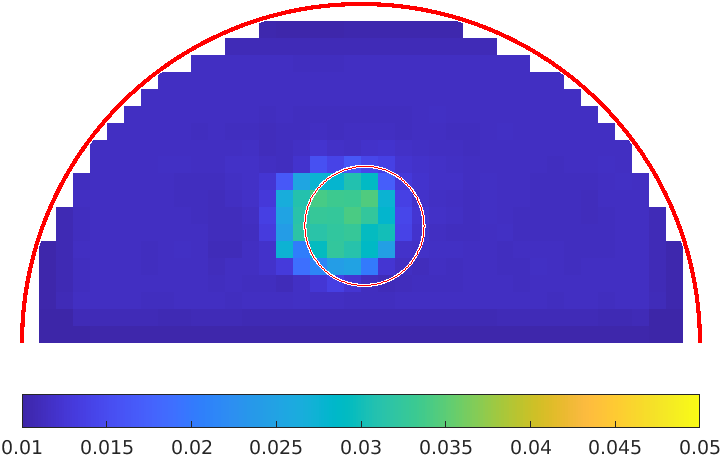}\hspace{0.025\textwidth}\includegraphics[width=\factor\textwidth]{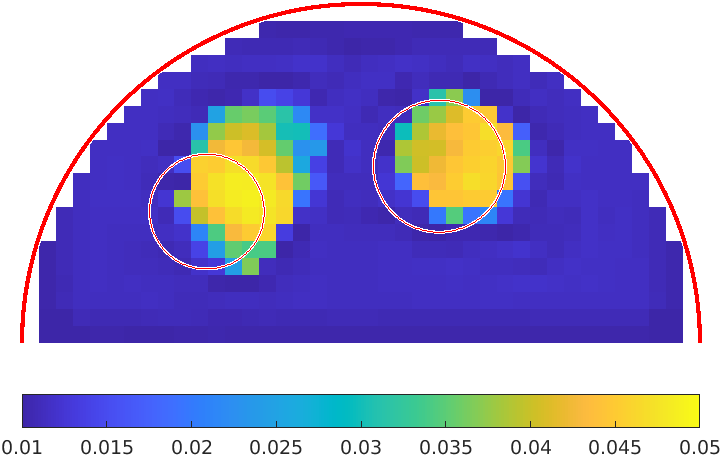}\hspace{0.025\textwidth}\includegraphics[width=\factor\textwidth]{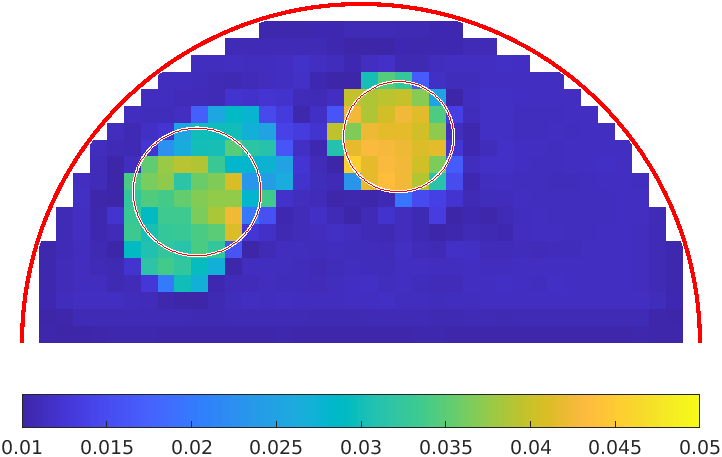}\hspace{0.025\textwidth}\includegraphics[width=\factor\textwidth]{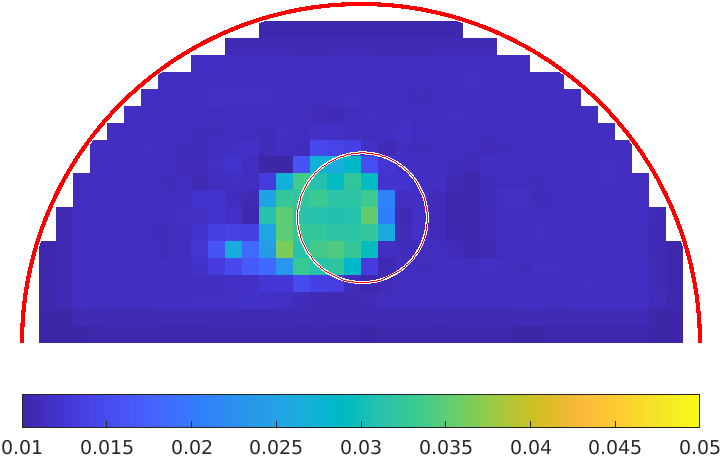}}
\end{center}
\caption{Reconstruction of the absorption coefficient distribution obtained with the Learned-SVD approach. \eqref{fig:gt} Ground truth. The 
red circles denote the true contrast regions along with the corresponding absorption coefficient, expressed 
as a multiplier of the background value. \eqref{fig:nfree}: noise-free case, \eqref{fig:G1}: 1\% noise, \eqref{fig:G3}: 3\% noise, \eqref{fig:G5}: 5\% noise.}
\label{fig:fResults}
\end{figure*}

\subsection{Learned-SVD: implementation of the network}
\noindent Both the encoders in the dAE and in sAE nets were chosen as fully connected layers coupled with a $\tanh$ activation function, while both the corresponding decoders 
were chosen as a  fully connected layer plus a sigmoid function. We investigated several options for 
the architecture of the $\Sigma$ bridge operator. Experiments showed that the best architecture was formed by 7 layers, each one consisting in a fully connected layer plus a $\tanh$ activation function.  
The training of the NNs was performed via the Stochastic Gradient Descent method, with default settings. The training strategy was the following: the dAE and the sAE 
were trained separately, with learning rate of 0.1 and 0.01, respectively. 
For each AE network, the loss function was chosen as the  Euclidean distance between the ground truth and the reconstructed data. The number of epochs was 2000 for each training procedure.
The $\Sigma$ network was trained on the dataset formed by the hidden representations of the two AE networks, 
with a learning rate 0.1. The final trained Learned-SVD 
network (cf. \Cref{fig:structures_BOINK}, green shaded area) was created by connecting the layers of the encoder of dAE, of the network $\Sigma$ and of the decoder of sAE.  
The Learned-SVD net provided good results, but artifacts were still present: in order to remove them, we added a downstream NN with the role of a further denoiser.  For this latter net we
used the same architecture of~\cite{Yoo2020}. An illustration of  the end--to--end workflow 
is provided in~\cref{fig:workflow}.

\subsection{Variational Approaches}
\label{ssec:varappr}
To carry out a comparison of the results obtained with the Learned-SVD approach with
classical variational methods we considered two approaches that we had used in our
past research work in the Rytov setting:

\paragraph{Elastic Net regularization~\cite{Causin19}} we look for 
\[
\argmin_{\mua} \frac12\|\vJ\,\mua-\vy\| + \alpha\lp\vartheta\|\mua\|_1+(1-\vartheta)\|\mua\|_2^2\rp,
\]
$\vartheta\in[0,1]$. In this case, $\mathcal{R}$ is a convex combination of $\ell_2$ and $\ell_1$ norms, which aims
at merging the best characteristics of both regularization functionals: the sparse--promoting and peak-enhancing property
of the $\ell_1$ norm and the robustness of the Tikhonov regularization. The solution is computed 
using the MatLab implementation of {\tt glmnet}~\cite{glmnet}, 
a highly efficient package that fits generalized linear and similar models via penalized maximum likelihood. 
Cross validation is used to obtain the optimal regularization parameter. In the numerical experiments $\vartheta$ is set to 0.5.

\paragraph{Bregman approach~\cite{Benfenati14,Benfenati2020}} it is an iterative strategy where at each step $\mathcal{R}$ is substituted with its Bregman Divergence computed at the previous iterate. If $\mathcal{R}$ satisfies suitable hypotheses 
(namely, being a convex, proper, lower semicontinuous function) and denoting by $\partial \mathcal{R}(\boldsymbol{\xi})$ 
the subdifferential of $\mathcal{R}$ computed at ${\boldsymbol{\xi}}$ (see \cite{rockafellar2015convex} for technical details), then the Bregman divergence of $\mathcal{R}$ is
\[
D_\mathcal{R}^{\vp}({\boldsymbol{\xi}}, {\boldsymbol{\eta}}) = \mathcal{R}({\boldsymbol{\xi}}) - \mathcal{R}({\boldsymbol{\eta}}) - \langle \vp,
{\boldsymbol{\xi}}-{\boldsymbol{\eta}}\rangle, \quad \vp \in\partial \mathcal{R}({\boldsymbol{\xi}}).
\]
The functional to minimize at the $k+1$--th step is thus
\begin{equation*}
\frac12\|\vJ\,\mua-\vy\|_2^2 + \alpha D_\mathcal{R}^{\vp^k}(\mua,\mua^k)
\end{equation*}
where $\vp^k\in\partial\mathcal{R}(\mua^k)$. Discarding the constant terms and considering the optimality conditions, the whole procedure consists in a two--step strategy
\begin{align}
\mua^{k+1} = \argmin_{\mua} \frac12\|\vJ\,\mua-\vy\|_2^2 +& \alpha \mathcal{R}(\mua)+\nonumber\label{eq:inner}\\
 &- \alpha \langle \vp^k, \mua\rangle\\
 \vp^{k+1} = \vp^k - \frac1\alpha \vJ^\top\lp\vJ\,\mua -\vy\rp
\end{align}
The Bregman procedure was coupled with the choice of $\ell_1$ regularization. The inner sub problem \eqref{eq:inner} is dealt
with the Forward--Backward algorithm \cite{9390397}, which consists in solving the general problem
$$
\min_{\boldsymbol{\eta}} f_0({\boldsymbol{\eta}}) + f_1({\boldsymbol{\eta}})
$$
via an iterative scheme
$$
\vx^{k+1} \gets \argmin_{\vx} f_0(\vx^k) + \frac{1}{2\gamma}\|\vx - \lp\vx^k - \gamma \nabla f_1(\vx)\rp\|_2^2  
$$
In view of \eqref{eq:inner}, $f_1(\cdot) = \frac12\|\vJ \cdot - y\|_2^2 -\alpha \la\cdot, \vp\ra$ and $f_0(\cdot) = \|\cdot\|_1$. The complete procedure is depicted in \cref{al:BregmanFB}, where $\soft({\boldsymbol{\eta}})$ 
stands for the component--wise soft thresholding $\soft_\beta(\eta) = \sign(\eta)
\max\{0,|{\eta}|-\beta\}$. The convergence of this approach has been proved in several previous works~\cite{doi:10.1137/070703983, doi:10.1137/090746379, Benfenati14,doi:10.1137/050626090}.

\begin{algorithm}
	\caption{Bregman Procedure with $\ell_1$ regularization}\label{al:BregmanFB}
	\begin{algorithmic}
		\STATE{Set $\mua^0$ such that $\vp^0=\bm{0}, \vp^0\in\partial \mathcal{R}(\mua^0)$, set $\alpha>0, \, \gamma>0$.}
		\FOR{$k=0,1,\dots$}
			\STATE{$\vx^\star \gets \mua^k$}
			\FOR{$l=0,1,\dots$}
				\STATE{$\vz^l \gets (\Id-\gamma\vJ^\top\vJ)\vx^\star + \gamma(\vJ^\top\vy + \vp^k)$}
				\STATE{$\vx^\star \gets \soft_{\gamma\alpha}(\vz^l)$}
			\ENDFOR
			\STATE{$\mua^k \gets \vx^\star$}
			\STATE{$\vp^{k+1}\gets \vp^k- \displaystyle\frac1\alpha \vJ^\top\lp\vJ\,\mua^k -\vy\rp$}
		\ENDFOR
	\end{algorithmic}
\end{algorithm}

The algorithm is stopped when the number or outer iterations reaches 100, whilst the inner solver runs for 50 iterations. The regularization parameter $\alpha$ was set to $1.5\,\|\vJ^\top\,\vy\|_\infty$ and $\gamma = 0.99/\|\vJ^\top\,\vJ\|$.

\begin{remark}
	With an abuse of notation, within the variational framework $\vy$ denotes the algorithm fluctuation $\log(u/u_0)$, recalling that ${u_0}$ are the measurements corresponding to unperturbed background conditions, and $u$ refers to condition (\emph{e.g.} compression of the tissue) in which the lightpath is perturbed by the presence of pathological regions.
\end{remark}

\subsection{Numerical results: qualitative and quantitative evaluation}
\begin{figure}[!t]
	\newcommand{\factor}{0.13}
	\begin{center}
		\subfloat[\label{fig:ela0}]{\includegraphics[width=\factor\textwidth]{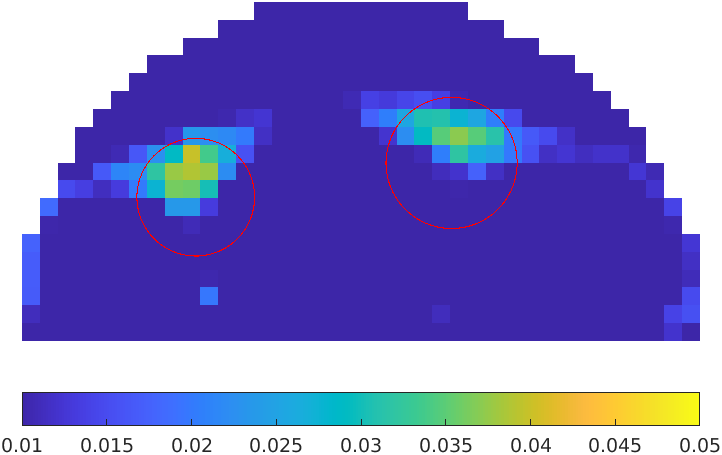}\hspace{0.025\textwidth}\includegraphics[width=\factor\textwidth]{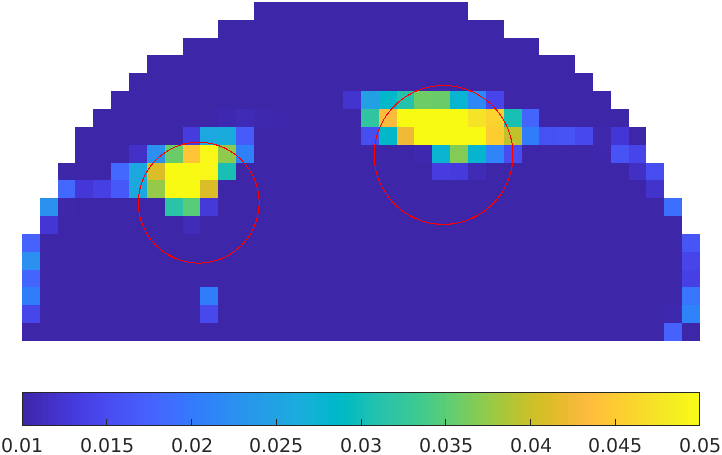}\hspace{0.025\textwidth}\includegraphics[width=\factor\textwidth]{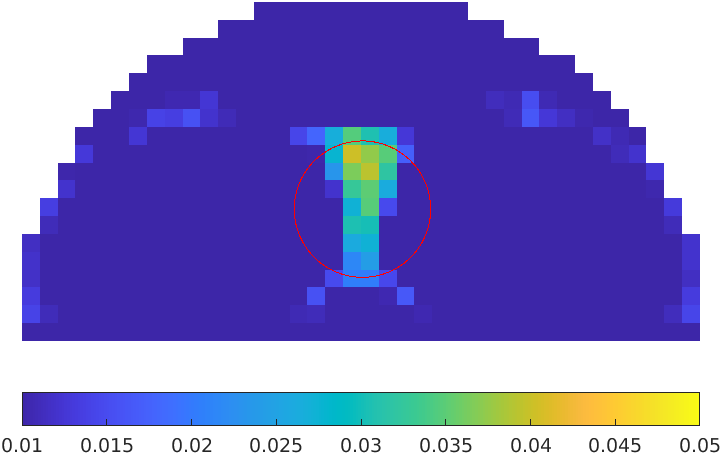}}
		
	\subfloat[\label{fig:ela1}]{\includegraphics[width=\factor\textwidth]{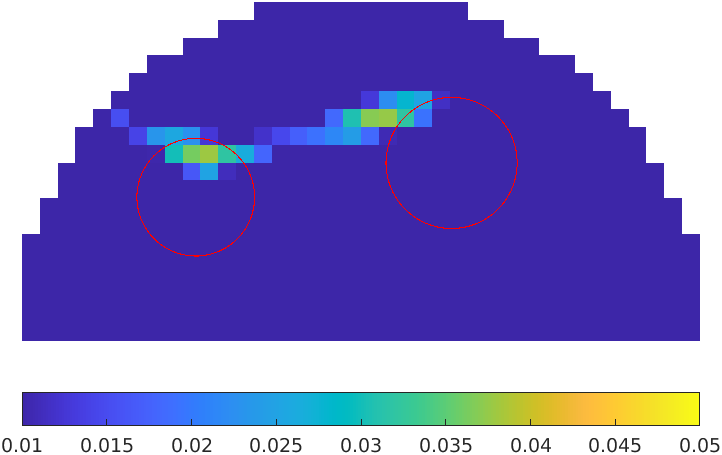}\hspace{0.025\textwidth}\includegraphics[width=\factor\textwidth]{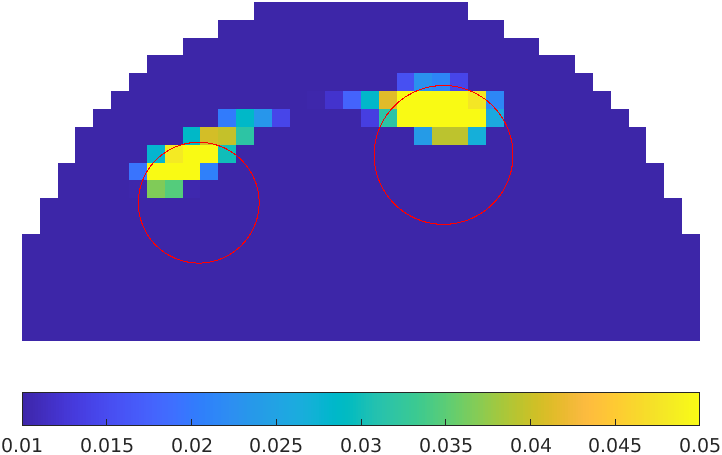}\hspace{0.025\textwidth}\includegraphics[width=\factor\textwidth]{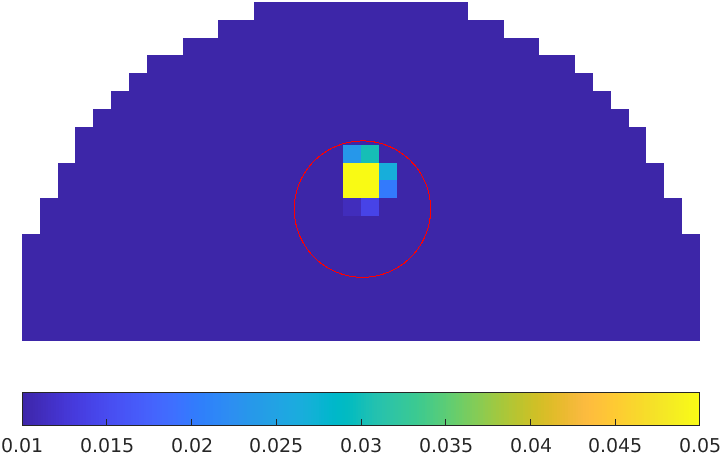}}

	\subfloat[\label{fig:ela3}]{\includegraphics[width=\factor\textwidth]{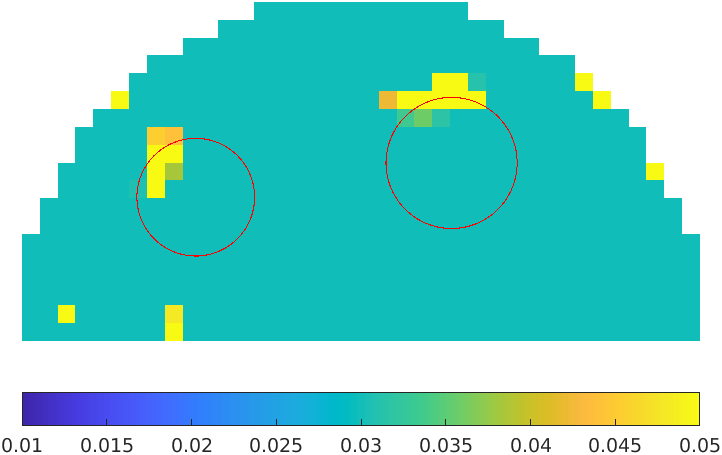}\hspace{0.025\textwidth}\includegraphics[width=\factor\textwidth]{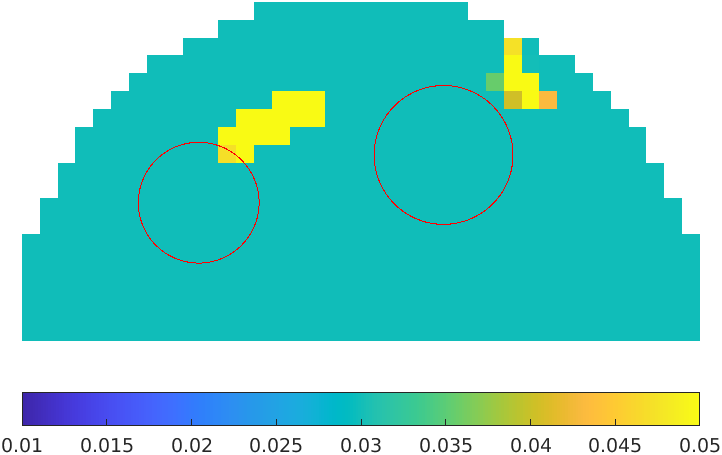}\hspace{0.025\textwidth}\includegraphics[width=\factor\textwidth]{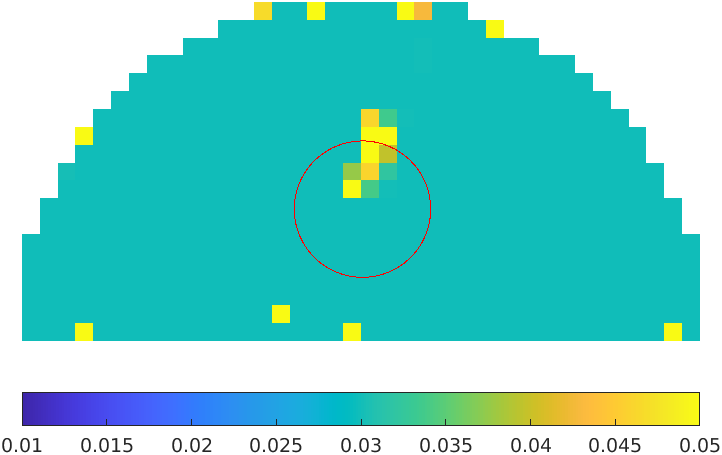}}

	\end{center}
	\caption{Reconstruction via the Elastic Net approach. \ref{fig:ela0}): noise-free case, \ref{fig:ela1}): 1\% noise, \ref{fig:ela3}): 3\% noise. The reconstruction are not reliable even in the noise free case due to the low resolution employed in these test.}
	\label{fig:ela}
\end{figure}

\begin{figure}[!t]
	\newcommand{\factor}{0.13}
	\begin{center}
		\subfloat[\label{fig:breg0}]{\includegraphics[width=\factor\textwidth]{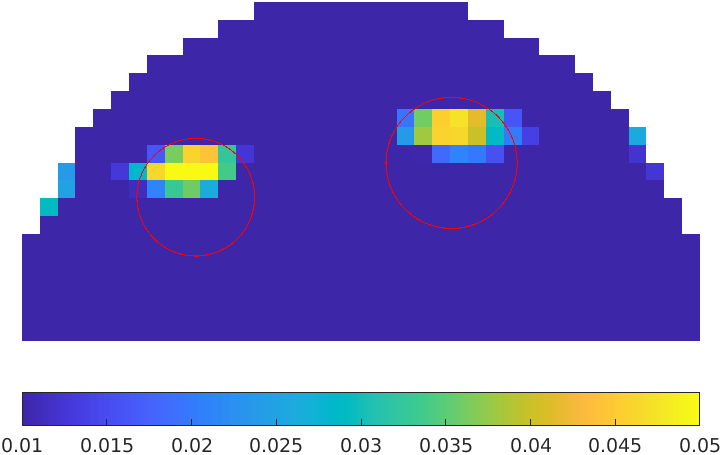}\hspace{0.025\textwidth}\includegraphics[width=\factor\textwidth]{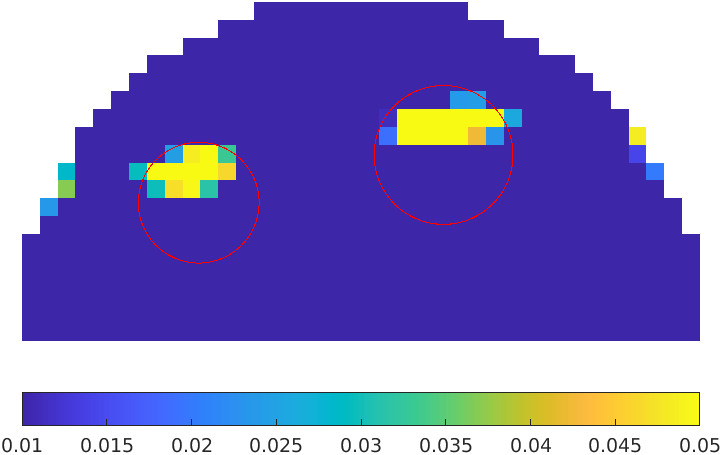}\hspace{0.025\textwidth}\includegraphics[width=\factor\textwidth]{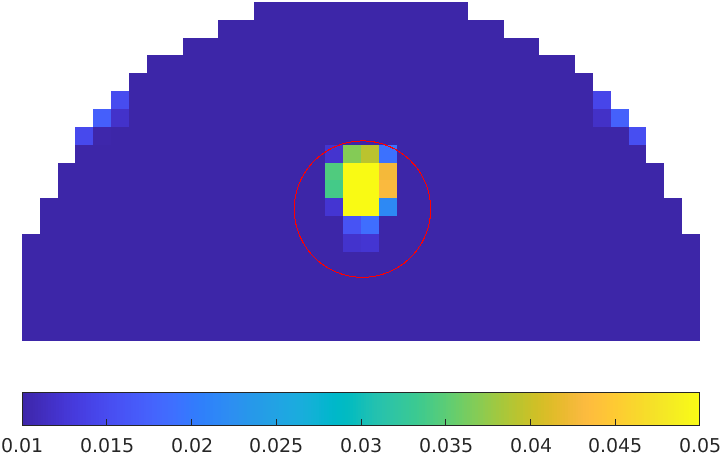}}
		
		\subfloat[\label{fig:breg1}]{\includegraphics[width=\factor\textwidth]{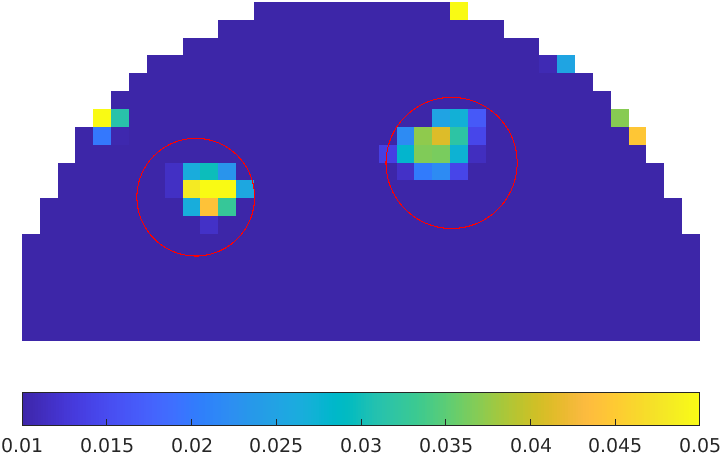}\hspace{0.025\textwidth}\includegraphics[width=\factor\textwidth]{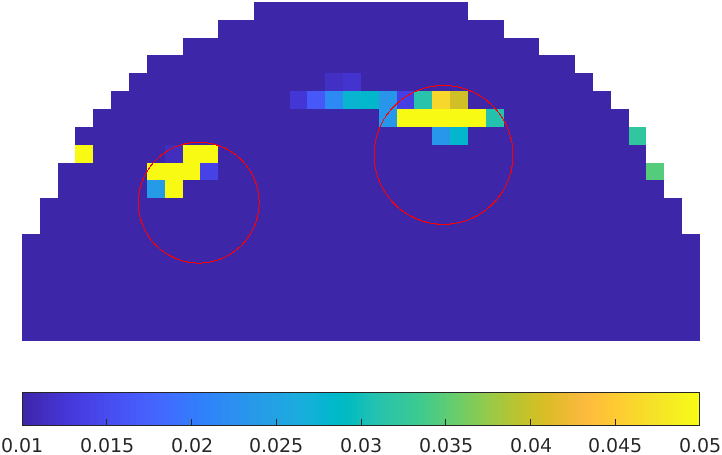}\hspace{0.025\textwidth}\includegraphics[width=\factor\textwidth]{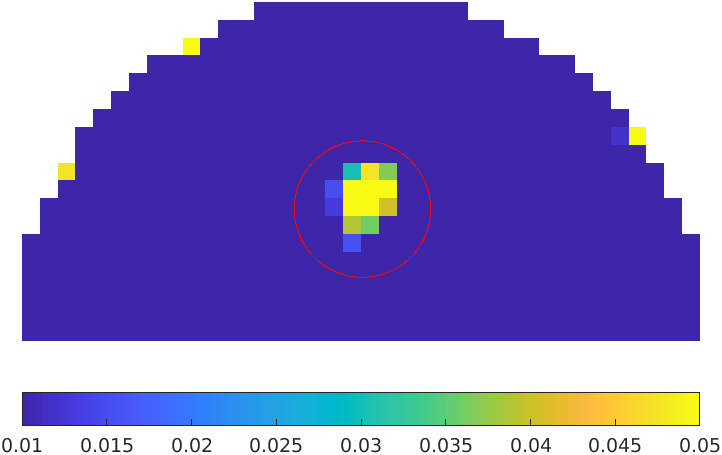}}
		
		\subfloat[\label{fig:breg3}]{\includegraphics[width=\factor\textwidth]{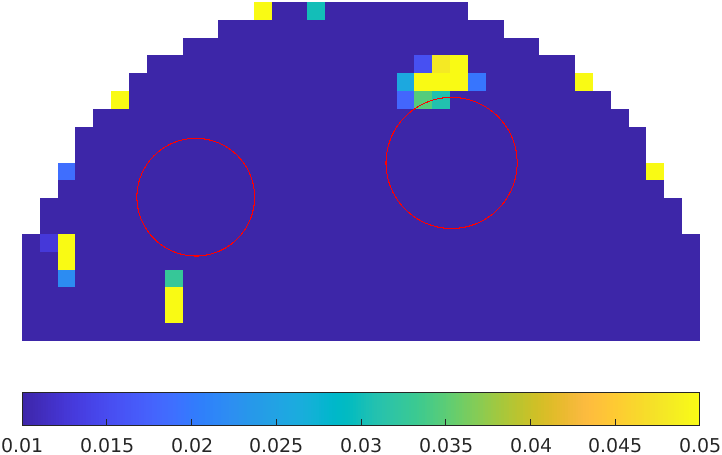}\hspace{0.025\textwidth}\includegraphics[width=\factor\textwidth]{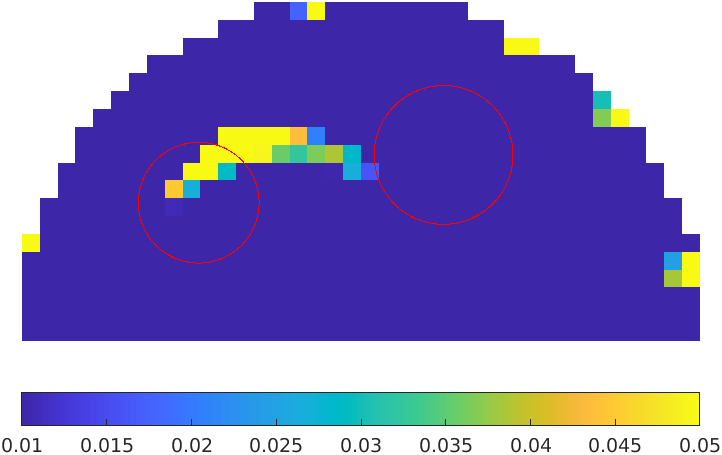}\hspace{0.025\textwidth}\includegraphics[width=\factor\textwidth]{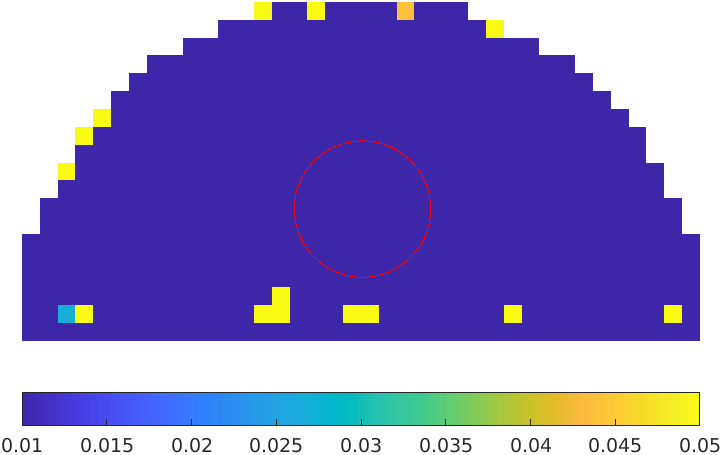}}
		
	\end{center}
	\caption{Reconstruction via the Bregman approach. \ref{fig:breg0}): noise-free case, 
	\ref{fig:breg1}): 1\% noise, \ref{fig:breg3}): 3\% noise.}
	\label{fig:breg}
\end{figure}

\noindent \Cref{fig:fResults} presents the results obtained for five different samples in the test set with varying levels of noise. The red circles depict the exact position and radius of the contrast region (ground--truth). In the noise-free case, the perturbed regions is almost perfectly recovered, both in location, shape and intensity. The other rows refer to the reconstruction obtained when Gaussian noise affects the  measure.  As one expects, the increasing noise level worsens the quality of the reconstruction but the gross location and shape of the perturbed region are correctly captured. 
The results obtained via the Learned-SVD approach are compared with the ones obtained by the 
two considered variational frameworks in \cref{fig:ela} (Elastic-Net) and in \cref{fig:breg} (Bregman), respectively. 
These figures show the results for test cases of columns 1, 3 and 5 in \cref{fig:gt}. Observe as 
already in the noise--free case (first row) we do not obtain truly reliable results due to the low resolution 
of the voxelization. Moreover, as the noise level increases the reconstruction quality 
significantly deteriorates, so that noise level 5\% was not reported due to complete lack of significance.

To quantitatively assess the performance of the proposed approach, we use two metrics. 
The first one is the average value of the absorption coefficient inside the reconstructed contrast regions (ACR),
which aims to assess the accuracy of the reconstructed amplitude of the contrast. 
To evaluate the ACR we first need to determine the extend of the computed contrast regions by recognizing the  
connected components in the reconstructed image and assigning each of them to the 
contrast region which minimizes the distance between with the centroid of the pixel cloud.
The second metric we use is the True Positive Ratio (TPR), which assesses the spatial
accuracy of the positioning of the reconstructed contrast regions 
by checking how many pixels in the reconstruction actually belong to the true 
contrast region. The results of the two metrics are summarized in~\Cref{tab:results}, which presents their averages 
on the 150 samples of the test set for the Learned-SVD approach and, for comparison,
of the considered variational approaches. Notice that the ACR results have been binned according to
the intensity of the contrast region. 
\begin{table*}[!t]
	\caption{ACR (binned according to
the intensity of the contrast region) and TPR metrics for different noise levels. Values
		are averaged over 150 test samples.}\label{tab:results}
\begin{center}
	\begin{tabular}{c|c|c|c|c}
		Noise Level & ${3\mu_{a,0}}$ (GT: 3e-2)&  ${4\mu_{a,0}}$ (GT: 4e-2) &  ${5\mu_{a,0}}$ (GT: 5e-02) & TPR\\
		\toprule
		&\multicolumn{4}{c}{Learned-SVD}\\
		\midrule
		0\% & 3.05e-02 $\pm$ 2.00e-03 & 4.08e-02 $\pm$ 3.72e-03 & 4.74e-02 $\pm$ 2.19e-03 & 0.94\\

		1\% & 3.07e-02 $\pm$ 2.77e-03 & 3.76e-02 $\pm$ 5.54e-03 & 4.38e-02 $\pm$ 4.01e-03 & 0.76\\

		3\% & 3.06e-02 $\pm$ 5.56e-03 & 3.61e-02 $\pm$ 6.10e-03 & 4.11e-02 $\pm$ 4.65e-03 & 0.52\\

		5\% & 3.00e-02 $\pm$ 3.82e-03 & 3.15e-02 $\pm$ 4.44e-03 & 3.39e-02 $\pm$ 6.04e-03 & 0.46\\
		
		\midrule		
		&\multicolumn{4}{c}{Elastic Net}\\
		\midrule
		0\% & 2.73e-02 $\pm$ 4.74e-03 & 3.45e-02 $\pm$ 4.96e-03 & 3.90e-02 $\pm$ 8.69e-03 & 0.45\\
		1\% & 2.91e-02 $\pm$ 9.36e-03 & 4.30e-02 $\pm$ 8.70e-03 & 5.01e-02 $\pm$ 9.80e-03 & 0.17\\
		3\% & 9.55e-02 $\pm$ 1.92e-02 & 1.25e-01 $\pm$ 2.83e-02 & 1.23e-01 $\pm$ 3.36e-02 & 0.05\\
            \midrule
&\multicolumn{4}{c}{Bregman}\\
		\midrule
		0\% & 4.02e-02 $\pm$ 8.95e-03 & 5.93e-02 $\pm$ 1.61e-02 & 8.54e-02 $\pm$ 2.98e-02 & 0.26\\
		1\% & 4.77e-02 $\pm$ 1.81e-02 & 5.85e-02 $\pm$ 1.84e-02 & 8.34e-02 $\pm$ 2.85e-02 & 0.17\\
		3\% & 1.28e-01 $\pm$ 5.50e-02 & 1.45e-01 $\pm$ 9.06e-02 & 1.36e-01 $\pm$ 6.80e-02 & 0.03\\
		\midrule		
	\end{tabular}
\end{center}

\end{table*}

For all the considered approaches, increasing levels of noise significantly worsen the TPR, 
as well as the ACR. However, the quality loss exhibited by the Learned-SVD approach 
is significantly less pronounced, especially for the ACR metric. 
As a matter of fact, the reconstructed regions with the NN-based method
still show a notable contrast with respect to the background value, and 46\% of the
extension of the perturbed regions is correctly recovered on average even for the highest noise level (vs  
5\% and 3\% for the variational approaches, respectively).
One should notice that in the Elastic Net approach the reconstructed intensities constantly
approximate the nominal one by defect for low noise level: this is not fully surprising, since 
regularization comes at the cost of a certain degree of smearing/blurring of the solution. On the other hand, Bregman technique increases the contrast, as already observed in previous works. In presence of high level of noise, the results are completely unreliable.
Classic variational methods suffer from this problem to a higher degree
especially for a rather coarse voxelization as the one we considered here (results, not reported
here, for finer voxelizations, do confirm in any case this general trend).

\section{Conclusions}
\label{sec:concl}
DOT reconstruction is a severely ill-conditioned problem which demands
a robust regularization to obtain reasonable results. 
Commonly used  strategies are based on physics-driven models
accompanied by priors which enforce constraints on the variance of the 
solution or in its sparsity.  These strategies often fail or are computationally
very intesive.  
In this work we have investigated
a NN-based regularization technique,
inspired by the  Learned-SVD method originally proposed in~\cite{Boink19}
for general inverse problem. 
This approach is a fully data--driven strategy which adopts two AEs 
bridged by an operator which mimics the effect of a (truncated) singular
value matrix. 
The Learned-SVD  produces significantly improved results with
respect to classic variational approaches based on $L^p, p \ge 1,$ penalization. 
According to the measured metrics, 
the present NN-based strategy yields an average 94\% TPR 
for noise--free data  and
for higher levels of noise the deterioration of the results
is signifiantly less severe than for variational methods, even for a coarse voxelization
of the domain.  
The training of the net performed in this work was 
based on a purely synthetic dataset: results obtained by other authors (see for 
example~\cite{sabir2020convolutional}) show 
that such an approach can convey good performances even when  
networks trained in such a way are used in realistic conditions. 
Transfer learning techniques could also be employed to
improve the reliability of the net on realistic datasets, even in presence of a limited
number of available samples. 
Future work will be focused on the idea of exporting and adapting techniques
from low-dose X-ray CT (another ill-conditioned problem) to the present context,
with  a specific attention for the hybrid integration with information coming from the rich 
physical models describing light propagation in tissues.  

\section{Acknowledgments.} AB and PC received support by the SEED-PRECISION project, funded by University of Milan. PC also acknowledges support from Italian Ministry of Research PRIN project NA\_FROM-PDEs.

\bibliographystyle{IEEEtran}
\bibliography{biblio}

\begin{thebibliography}{10}
\providecommand{\url}[1]{#1}
\csname url@samestyle\endcsname
\providecommand{\newblock}{\relax}
\providecommand{\bibinfo}[2]{#2}
\providecommand{\BIBentrySTDinterwordspacing}{\spaceskip=0pt\relax}
\providecommand{\BIBentryALTinterwordstretchfactor}{4}
\providecommand{\BIBentryALTinterwordspacing}{\spaceskip=\fontdimen2\font plus
\BIBentryALTinterwordstretchfactor\fontdimen3\font minus
  \fontdimen4\font\relax}
\providecommand{\BIBforeignlanguage}[2]{{%
\expandafter\ifx\csname l@#1\endcsname\relax
\typeout{** WARNING: IEEEtran.bst: No hyphenation pattern has been}%
\typeout{** loaded for the language `#1'. Using the pattern for}%
\typeout{** the default language instead.}%
\else
\language=\csname l@#1\endcsname
\fi
#2}}
\providecommand{\BIBdecl}{\relax}
\BIBdecl

\bibitem{hoshi2016overview}
Y.~Hoshi and Y.~Yamada, ``Overview of diffuse optical tomography and its
  clinical applications,'' \emph{Journal of biomedical optics}, vol.~21, no.~9,
  p. 091312, 2016.

\bibitem{yamada2014diffuse}
Y.~Yamada and S.~Okawa, ``Diffuse optical tomography: Present status and its
  future,'' \emph{Optical Review}, vol.~21, no.~3, pp. 185--205, 2014.

\bibitem{jiang2018diffuse}
H.~Jiang, \emph{Diffuse Optical Tomography: Principles and Applications}.\hskip
  1em plus 0.5em minus 0.4em\relax CRC Press, 2018.

\bibitem{arridge1999optical}
S.~R. Arridge, ``Optical tomography in medical imaging,'' \emph{Inverse
  problems}, vol.~15, no.~2, p. R41, 1999.

\bibitem{Boas01}
D.~Boas, D.~Brooks, E.~Miller, C.~DiMarzio, M.~Kilmer, R.~Gaudette, and
  Q.~Zhang, ``Imaging the body with diffuse optical tomography,'' \emph{IEEE
  Signal Processing Magazine}, vol.~18, no.~6, pp. 57--75, 2001.

\bibitem{Arridge2009optical}
S.~R. Arridge and J.~C. Schotland, ``Optical tomography: forward and inverse
  problems,'' \emph{Inverse problems}, vol.~25, no.~12, p. 123010, 2009.

\bibitem{5728925}
O.~Lee, J.~M. Kim, Y.~Bresler, and J.~C. Ye, ``Compressive diffuse optical
  tomography: Noniterative exact reconstruction using joint sparsity,''
  \emph{IEEE Transactions on Medical Imaging}, vol.~30, no.~5, pp. 1129--1142,
  2011.

\bibitem{doi:10.1137/090781590}
H.~Egger and M.~Schlottbom, ``Analysis and regularization of problems in
  diffuse optical tomography,'' \emph{SIAM Journal on Mathematical Analysis},
  vol.~42, no.~5, pp. 1934--1948, 2010.

\bibitem{10.1117/12.2230074}
\BIBentryALTinterwordspacing
H.~O. Kazanci and S.~L. Jacques, ``{Diffuse light tomography to detect blood
  vessels using Tikhonov regularization},'' in \emph{Saratov Fall Meeting 2015:
  Third International Symposium on Optics and Biophotonics and Seventh
  Finnish-Russian Photonics and Laser Symposium (PALS)}, E.~A. Genina, V.~V.
  Tuchin, V.~L. Derbov, D.~E. Postnov, I.~V. Meglinski, K.~V. Larin, and A.~B.
  Pravdin, Eds., vol. 9917, International Society for Optics and
  Photonics.\hskip 1em plus 0.5em minus 0.4em\relax SPIE, 2016, pp. 202 -- 210.
  [Online]. Available: \url{https://doi.org/10.1117/12.2230074}
\BIBentrySTDinterwordspacing

\bibitem{Okawa:11}
\BIBentryALTinterwordspacing
S.~Okawa, Y.~Hoshi, and Y.~Yamada, ``Improvement of image quality of
  time-domain diffuse optical tomography with lp sparsity regularization,''
  \emph{Biomed. Opt. Express}, vol.~2, no.~12, pp. 3334--3348, Dec 2011.
  [Online]. Available:
  \url{http://opg.optica.org/boe/abstract.cfm?URI=boe-2-12-3334}
\BIBentrySTDinterwordspacing

\bibitem{Causin19}
P.~Causin, G.~Naldi, and R.~Weishaeupl, ``Elastic net regularization in diffuse
  optical tomography applications,'' in \emph{2019 IEEE 16th International
  Symposium on Biomedical Imaging (ISBI 2019)}, 2019, pp. 1627--1630.

\bibitem{Causin20}
P.~Causin, M.~G. Lupieri, G.~Naldi, and R.-M. Weishaeupl, ``Mathematical and
  numerical challenges in optical screening of female breast,''
  \emph{International Journal for Numerical Methods in Biomedical Engineering},
  vol.~36, no.~2, p. e3286, 2020.

\bibitem{Benfenati2020}
A.~Benfenati, P.~Causin, M.~Lupieri, and G.~Naldi, ``Regularization techniques
  for inverse problem in {DOT} applications,'' \emph{Journal of Physics:
  Conference Series}, vol. 1476, p. 012007, mar 2020.

\bibitem{Zhang201120}
X.~Zhang, M.~Burger, and S.~Osher, ``A unified primal-dual algorithm framework
  based on {B}regman iteration,'' \emph{Journal of Scientific Computing},
  vol.~46, no.~1, pp. 20--46, 2011.

\bibitem{Benfenati14}
A.~Benfenati and V.~Ruggiero, ``Inexact {B}regman iteration with an application
  to {P}oisson data reconstruction,'' \emph{Inverse Problems}, vol.~29, no.~6,
  2013.

\bibitem{Benfenati6}
A.~Benfenati, A.~La~Camera, and M.~Carbillet, ``Deconvolution of post-adaptive
  optics images of faint circumstellar environments by means of the inexact
  {B}regman procedure,'' \emph{Astronomy and Astrophysics}, vol. 586, 2016.

\bibitem{Benfenati11}
A.~Benfenati and V.~Ruggiero, ``Inexact {B}regman iteration for deconvolution
  of superimposed extended and point sources,'' \emph{Communications in
  Nonlinear Science and Numerical Simulation}, vol.~20, no.~3, pp. 882--896,
  2015.

\bibitem{Schweiger2014toast++}
M.~Schweiger and S.~R. Arridge, ``The {TOAST}++ software suite for forward and
  inverse modeling in optical tomography,'' \emph{J Biomed Opt}, vol.~19,
  no.~4, p. 040801, 2014.

\bibitem{Wu20}
Z.~Wu, Y.~Sun, A.~Matlock, J.~Liu, L.~Tian, and U.~S. Kamilov, ``{SIMBA}:
  Scalable inversion in optical tomography using deep denoising priors,''
  \emph{IEEE Journal of Selected Topics in Signal Processing}, vol.~14, no.~6,
  pp. 1163--1175, 2020.

\bibitem{Yoo2020}
J.~Yoo, S.~Sabir, D.~Heo, K.~Kim, A.~Wahab, Y.~Choi, S.-I. Lee, E.~Chae,
  H.~Kim, Y.~Bae, Y.-W. Choi, S.~Cho, and J.~Ye, ``Deep learning diffuse
  optical tomography,'' \emph{IEEE Transactions on Medical Imaging}, vol.~39,
  no.~4, pp. 877--887, 2020.

\bibitem{mozumder2021model}
M.~Mozumder, A.~Hauptmann, I.~Nissil{\"a}, S.~R. Arridge, and T.~Tarvainen, ``A
  model-based iterative learning approach for diffuse optical tomography,''
  \emph{arXiv preprint arXiv:2104.09579}, 2021.

\bibitem{zhang2019brief}
L.~Zhang and G.~Zhang, ``Brief review on learning-based methods for optical
  tomography,'' \emph{Journal of Innovative Optical Health Sciences}, vol.~12,
  no.~06, p. 1930011, 2019.

\bibitem{Boink19}
Y.~E. Boink and C.~Brune, ``Learned {SVD:} solving inverse problems via hybrid
  autoencoding,'' \emph{CoRR}, vol. abs/1912.10840, 2019.

\bibitem{Lunz18}
S.~Lunz, O.~Öktem, and C.-B. Schönlieb, ``Adversarial regularizers in inverse
  problems,'' vol. 2018-December, 2018, p. 8507 – 8516.

\bibitem{Kobler_2020_CVPR}
E.~Kobler, A.~Effland, K.~Kunisch, and T.~Pock, ``Total deep variation for
  linear inverse problems,'' in \emph{Proceedings of the IEEE/CVF Conference on
  Computer Vision and Pattern Recognition (CVPR)}, June 2020.

\bibitem{Kobler20}
------, ``Total deep variation: {A} stable regularizer for inverse problems,''
  \emph{CoRR}, vol. abs/2006.08789, 2020.

\bibitem{Li2020}
H.~Li, J.~Schwab, S.~Antholzer, and M.~Haltmeier, ``{NETT}: solving inverse
  problems with deep neural networks,'' \emph{Inverse Problems}, vol.~36,
  no.~6, p. 065005, jun 2020.

\bibitem{vavadi2016automated}
H.~Vavadi and Q.~Zhu, ``Automated data selection method to improve robustness
  of diffuse optical tomography for breast cancer imaging,'' \emph{Biomedical
  optics express}, vol.~7, no.~10, pp. 4007--4020, 2016.

\bibitem{Durduran2010DiffuseOF}
T.~Durduran, R.~Choe, W.~B. Baker, and A.~G. Yodh, ``Diffuse optics for tissue
  monitoring and tomography.'' \emph{Reports on progress in physics. Physical
  Society}, vol. 73 7, 2010.

\bibitem{Mansuripur2002}
M.~Mansuripur, \emph{Classical optics and its applications}.\hskip 1em plus
  0.5em minus 0.4em\relax Cambridge University Press, 2002.

\bibitem{Ishimaru78}
A.~Ishimaru, \emph{Wave propagation and scattering in random media}.\hskip 1em
  plus 0.5em minus 0.4em\relax Academic Press New York, 1978, vol.~2.

\bibitem{causin2020inverse}
P.~Causin and R.-M. Weishaeupl, ``Inverse problems in diffuse optical
  tomography applications,'' in \emph{Mathematical Modelling in Real Life
  Problems}.\hskip 1em plus 0.5em minus 0.4em\relax Springer, 2020, pp. 1--16.

\bibitem{Sciacca21}
G.~Di~Sciacca, L.~Di~Sieno, A.~Farina, P.~Lanka, E.~Venturini, P.~Panizza,
  A.~Dalla~Mora, A.~Pifferi, P.~Taroni, and S.~Arridge, ``Enhanced diffuse
  optical tomographic reconstruction using concurrent ultrasound information,''
  \emph{Philosophical Transactions of the Royal Society A}, vol. 379, no. 2204,
  p. 20200195, 2021.

\bibitem{Panagiotou:09}
C.~Panagiotou, S.~Somayajula, A.~P. Gibson, M.~Schweiger, R.~M. Leahy, and
  S.~R. Arridge, ``Information theoretic regularization in diffuse optical
  tomography,'' \emph{J. Opt. Soc. Am. A}, vol.~26, no.~5, pp. 1277--1290, May
  2009.

\bibitem{Cao:07}
\BIBentryALTinterwordspacing
N.~Cao, A.~Nehorai, and M.~Jacob, ``Image reconstruction for diffuse optical
  tomography using sparsity regularization and expectation-maximization
  algorithm,'' \emph{Opt. Express}, vol.~15, no.~21, pp. 13\,695--13\,708, Oct
  2007. [Online]. Available:
  \url{http://opg.optica.org/oe/abstract.cfm?URI=oe-15-21-13695}
\BIBentrySTDinterwordspacing

\bibitem{meinshausen2006high}
N.~Meinshausen and P.~B{\"u}hlmann, ``High-dimensional graphs and variable
  selection with the lasso,'' \emph{The annals of statistics}, vol.~34, no.~3,
  pp. 1436--1462, 2006.

\bibitem{yodh1995spectroscopy}
A.~Yodh and B.~Chance, ``Spectroscopy and imaging wih diffusing light,''
  \emph{Physics Today}, vol.~48, no.~3, pp. 34--40, 1995.

\bibitem{glmnet}
\BIBentryALTinterwordspacing
J.~Friedman, T.~Hastie, and R.~Tibshirani, ``Regularization paths for
  generalized linear models via coordinate descent,'' \emph{Journal of
  Statistical Software}, vol.~33, no.~1, pp. 1--22, 2010. [Online]. Available:
  \url{https://www.jstatsoft.org/v33/i01/}
\BIBentrySTDinterwordspacing

\bibitem{rockafellar2015convex}
R.~Rockafellar, \emph{Convex Analysis: (PMS-28)}, ser. Princeton Landmarks in
  Mathematics and Physics.\hskip 1em plus 0.5em minus 0.4em\relax Princeton
  University Press, 2015.

\bibitem{9390397}
P.~L. Combettes and J.-C. Pesquet, ``Fixed point strategies in data science,''
  \emph{IEEE Transactions on Signal Processing}, vol.~69, pp. 3878--3905, 2021.

\bibitem{doi:10.1137/070703983}
W.~Yin, S.~Osher, D.~Goldfarb, and J.~Darbon, ``Bregman iterative algorithms
  for $\ell\_1$-minimization with applications to compressed sensing,''
  \emph{SIAM Journal on Imaging Sciences}, vol.~1, no.~1, pp. 143--168, 2008.

\bibitem{doi:10.1137/090746379}
X.~Zhang, M.~Burger, X.~Bresson, and S.~Osher, ``Bregmanized nonlocal
  regularization for deconvolution and sparse reconstruction,'' \emph{SIAM
  Journal on Imaging Sciences}, vol.~3, no.~3, pp. 253--276, 2010.

\bibitem{doi:10.1137/050626090}
P.~L. Combettes and V.~R. Wajs, ``Signal recovery by proximal forward-backward
  splitting,'' \emph{Multiscale Modeling \& Simulation}, vol.~4, no.~4, pp.
  1168--1200, 2005.

\bibitem{sabir2020convolutional}
S.~Sabir, S.~Cho, Y.~Kim, R.~Pua, D.~Heo, K.~H. Kim, Y.~Choi, and S.~Cho,
  ``Convolutional neural network-based approach to estimate bulk optical
  properties in diffuse optical tomography,'' \emph{Applied Optics}, vol.~59,
  no.~5, pp. 1461--1470, 2020.

\end{thebibliography}

\end{document}